\documentclass[10pt]{article}
\title{On Star Monotone Operators: Theory and Applications}
\author{Phan Quoc Khanh\thanks{Analytical and Algebraic Methods in Optimization Research Group (AAMO), Faculty of Mathematics and Statistics, Ton Duc Thang University, Ho Chi Minh City, Vietnam. Email: phanquockhanh@tdtu.edu.vn; pqkhanhus@gmail.com, ORCID-ID: 0000-0002-0524-0498}
\and
Felipe Lara\thanks{Instituto de Alta investigaci\'on (IAI), Universidad de Tarapac\'a, 
Arica, Chile. E-mail: felipelaraobreque@gmail.com; flarao@academicos.uta.cl. 
Web: felipelara.cl, ORCID-ID: 0000-0002-9965-0921} }

\usepackage{amsmath}
\usepackage{amsthm}
\usepackage{amssymb}
\usepackage{multicol}
\usepackage[active]{srcltx}
\usepackage{algorithm}
\usepackage{array}
\usepackage[dvipsnames,svgnames]{xcolor}
\usepackage{url}
\usepackage{tikz}
\usepackage{subfig}
\usepackage{graphicx}
\usepackage{pgfplots}
\usepackage{xcolor}
\usepackage{ulem}  
\usepackage{comment}
\providecommand{\U}[1]{\protect \rule{.1in}{.1in}}
\newtheorem{theorem}{Theorem}

\newtheorem{corollary}[theorem]{Corollary}

\newtheorem{definition}[theorem]{Definition}

\newtheorem{lemma}[theorem]{Lemma}

\newtheorem{proposition}[theorem]{Proposition}
\newtheorem{remark}[theorem]{Remark}

\newcommand{\R}{\mathbb{R}}

\newcommand{\zer}{\operatorname{zer}}

\newcommand{\ran}{\operatorname{ran}}

\begin{document}

%\today

\maketitle

%\centerline{\today}

\begin{abstract}
\noindent Introduced by Solodov and Svaiter \cite{SS}, star monotone operators have repeatedly appeared as a sufficient condition ensuring global convergence of first-order algorithms, without requiring full monotonicity. Yet, despite this recurring role, a systematic study of the class itself has been missing from the literature. This paper fills that gap with the first comprehensive treatment of star monotone operators, a class that strictly encompasses monotone and pseudomonotone operators with nonempty zero sets, together with operators that are strongly quasimonotone along rays and satisfy a radial nondegeneracy condition. We develop basic calculus rules, existence and uniqueness results for zeros, global metric subregularity, and precise structural relationships with classical monotonicity notions. We further demonstrate the versatility of star monotonicity through a Clarke subdifferential characterization of nonsmooth star quasiconvexity, applications to Minty and Stampacchia variational inequalities, linear convergence guarantees for subgradient-type methods, and separable operators related to the Debreu--Koopmans problem.

{\small}

\medskip


\medskip

\end{abstract}

\section{Introduction}

The study of nonlinear operator equations $0 \in T(x)$ occupies a central place in nonlinear functional analysis, tracing back to the foundational existence and surjectivity results for monotone operators on Hilbert and Banach spaces developed by Zarantonello, Minty, Browder and Rockafellar \cite{Browder1965, Minty1962, Rocka-1970, Zaran}, whose modern operator-theoretic treatment is presented systematically in \cite{BC}. This body of work, together with further studies for pseudomonotone and quasimonotone operators \cite{HKS,HSW}, underlies the modern treatment of variational inequalities and remains a driving force in the analysis of nonlinear equations, inclusions, and their associated fixed-point and resolvent operators. Beyond their role within functional analysis itself, monotone operator methods and variational inequalities have proved indispensable across a wide range of applied fields: they originated in the study of free boundary and obstacle problems for elliptic partial differential equations \cite{LS,Minty1967}, and have since become a standard tool in the analysis of nonlinear elliptic and parabolic PDEs, contact mechanics, optimal control, and mathematical economics, where they formalize equilibrium and complementarity conditions.

Parallel to this operator-theoretic development, generalized convexity has played an analogous unifying role on the function side. Quasiconvexity (the property that every sublevel set is convex) generalizes convexity while preserving the tractability needed for existence and characterization results, and is the mathematical formalization of the {\it tendency to diversification} assumption in consumer preference theory in microeconomics \cite{D-1959}, remaining fundamental in production theory and general equilibrium analysis \cite{FOS,Mas}. The classical correspondence (convex functions have monotone gradients, quasiconvex functions have quasimonotone gradients, see \cite{HKS}) links these two strands of the theory directly.

Within this landscape, star quasiconvexity has recently emerged as a considerably weaker notion, anchored at a single minimizer rather than requiring a two-point comparison: a function $h$ is star quasiconvex with respect to a minimizer $\overline{x}$ when the defining inequality of quasiconvexity is only imposed between $\overline{x}$ and every other point, rather than between every pair of points \cite{KL}. This single relaxation is enough to recover sharp characterizations (star-shapedness of sublevel sets, monotonicity along rays, separability, and a first-order criterion) while strictly enlarging the class of functions covered, and it resolves the converse of the classical Debreu--Koopmans theorem on additively decomposed quasiconvex functions \cite{CL,DK,Lara-10}.

The natural operator-theoretic counterpart of this idea (imposing the defining inequality of monotonicity only between a fixed reference point $x_0$ and every other point, rather than between every pair) is what we call {\it star monotonicity}. Although the resulting notion is new as a systematically studied object, the underlying inequality is not: it has repeatedly surfaced, usually without a name, as a sufficient condition for global convergence of first-order methods for variational inequality problems, going back at least to Solodov and Svaiter \cite{SS}, and appearing since in a variety of algorithmic contexts \cite{CGL2023,KLL2022,Mal}. Despite this recurring role, no systematic operator-theoretic treatment of star monotone operators (their basic calculus, existence and uniqueness theory, structural properties, and precise relationship to the classical generalized monotonicity hierarchy) has previously been carried out. This paper fills that gap, and in doing so lays the groundwork for extending Zarantonello- and Minty-type surjectivity results to this broader class of operators.

\medskip

{\bf Our contribution:} We introduce star monotone operators as a systematically studied class and develop a comprehensive operator-theoretic treatment. We establish their basic calculus (scaling, translation, weighted sums, stability under perturbations), existence and uniqueness of zeros under weak hemicontinuity, a global metric subregularity estimate, and a Tikhonov-type regularization scheme. We give a precise structural analysis of star monotonicity relative to the classical monotone/pseudomonotone/quasimonotone hierarchy, including a sharp ray-wise characterization in terms of quasimonotonicity restricted to rays through the anchor point, and we identify genuine structural limitations of the class by showing that star monotone operators need not have convex, closed, or bounded values. We further demonstrate the versatility of the notion through four applications: a Clarke subdifferential characterization of nonsmooth (strongly) star quasiconvex functions \cite{clarke,FLG}; a complete equivalence between zeros, star monotonicity, and solutions to the Minty and Stampacchia variational inequalities, extending classical results \cite{CGR,Lara-7} to this setting under weak hemicontinuity at a single point; linear convergence guarantees for a fixed-point scheme and for the subgradient method applied to nonsmooth star quasiconvex minimization; and a separable characterization of star monotone operators that partially extends the classical results of Crouzeix and Hassouni \cite{CH1,CH2} on quasimonotonicity of separable products, in the spirit of the Debreu--Koopmans problem \cite{DK,Lara-10}.

\medskip

The paper is organized as follows. Section \ref{sec:02} recalls the necessary preliminaries on generalized convexity and generalized monotonicity, including star quasiconvexity \cite{KL,Lara-10,NTV} and the classical monotone/pseudomonotone/quasi\-mo\-no\-to\-ne hierarchy \cite{HKS,HSW}, together with the continuity notions for set-valued operators used throughout the paper. Section \ref{sec:03} introduces star monotone operators and develops their basic theory, culminating in the existence, uniqueness, and ray-wise characterization results described above. Section \ref{sec:04} presents the four applications summarized above: the Clarke subdifferential characterization, the Minty--Stampacchia equivalence, the linear convergence results, and the separable characterization related to the Debreu--Koopmans problem. Finally, in Section \ref{sec:05}, we describe conclusions and further research directions.

\section{Preliminaries and Basic Definitions} \label{sec:02}

The inner product in $\mathbb{R}^{n}$ and the Euclidean norm are denoted by $\langle \cdot,\cdot \rangle$ and $\lVert \cdot \rVert$, respectively. %We denote  $\mathbb{R}_{+} := [0, + \infty[$ and $ \mathbb{R}_{++} := \, ]0, + \infty[$. %and $\mathbb{N}_{0} := \mathbb{N} \cup \{0\}$. 
Let $K$ be a nonempty set in $\mathbb{R}^{n}$, its topological interior by ${\rm int}\,K$ and its convex hull by ${\rm co}\,K$. %and its closure by $\overline{K}$. 
The open ball with center at $x_{0}$ and radius $r>0$ is denoted by $\mathbb{B} (x_{0}, r)$.

Let $K \subset \mathbb{R}^{n}$ be nonempty and $x_{0} \in K$.  $K$ is said to be star-shaped at the point $x_{0}$ if for every $x \in K$, the segment $[x_{0}, x] \subset K$.

%Given any $x, y, z \in \mathbb{R}^{n}$ and any $\beta \in \mathbb{R}$, the following relations hold:
%\begin{align}
%  \lVert \beta x + (1-\beta) y \rVert^{2} &= \beta \lVert x \rVert^{2} + (1 - \beta) \lVert y\rVert^{2} - \beta(1 - \beta) \lVert x - y \rVert^{2}, \label{iden:1}\\
% \langle x - z, y - x \rangle &= \frac{1}{2} \lVert z - y \rVert^{2} - \frac{1}{2} \lVert x - z \rVert^{2} - \frac{1}{2} \lVert y - x \rVert^{2}. \label{3:points}
%\end{align}

Given any extended-valued function $h: \mathbb{R}^{n} \rightarrow \overline{\mathbb{R}} := \mathbb{R} \cup \{\pm \infty\}$, the effective domain of $h$ is defined by ${\rm dom}\,h := \{x \in \mathbb{R}^{n}: h(x) < + \infty \}$.  It is said that $h$ is proper if ${\rm dom}\,h$ is nonempty and $h(x) > - \infty$ for all $x \in \mathbb{R}^{n}$. The notion of properness is important when dealing with minimization pro\-blems.
The set ${\rm epi}\,h := \{(x,t) \in \mathbb{R}^{n} \times \mathbb{R}: h(x) \leq t\}$ is the epigraph of $h$, $S_{\delta} (h) := \{x \in \mathbb{R}^{n}: h(x) \leq \delta\}$ is the sublevel set of $h$ at the height $\delta \in \mathbb{R}$, and  ${\rm argmin}_{ \mathbb{R}^{n}} h$ is the set of all minimal points of $h$.

A function $h$ with convex domain is said to be:
\begin{itemize}
 \item[$(a)$] (strongly) convex on ${\rm dom}\,h$ when there exists $\gamma \geq 0$ such that if for all $x, y \in \mathrm{dom}\,h$ and all $\lambda \in [0, 1]$, we have
 \begin{equation*}\label{strong:convex}
  h(\lambda y + (1-\lambda)x) \leq \lambda h(y) + (1-\lambda) h(x) - \lambda (1 - \lambda) \frac{\gamma}{2} \lVert x - y \rVert^{2},
 \end{equation*}

 \item[$(b)$] (strongly) quasiconvex on ${\rm dom}\,h$ when there exists $\gamma \geq 0$ such that if for all $x, y \in \mathrm{dom}\,h$ and all $\lambda \in [0, 1]$, we have
 \begin{equation*}\label{strong:quasiconvex}
  h(\lambda y + (1-\lambda)x) \leq \max \{h(y), h(x)\} - \lambda(1 -
  \lambda) \frac{\gamma}{2} \lVert x - y \rVert^{2}.
 \end{equation*}
\end{itemize} 
\noindent Every (strongly) convex function is (strongly) quasiconvex, while the reverse statement does not hold (see \cite{CM-Book,HKS,Lara-9}). Furthermore, when we say that a function is {\it strongly convex (resp., quasiconvex) with modulus $\gamma \geq 0$} we refer to both strongly convex (resp., quasiconvex) when $\gamma > 0$ and convex (resp., qua\-siconvex) when $\gamma = 0$.

Let $h: \mathbb{R}^{n} \rightarrow \mathbb{R}$ be a differentiable function. Then $h$ is said to be pseudoconvex (see \cite{Manga}) if
\begin{equation}\label{pseudoconvex}
 \langle \nabla h(x), y - x \rangle \geq 0 ~ \Longrightarrow ~ h(y) \geq h(x).
\end{equation}
If $h$ is pseudoconvex, then every local minimum is global minimum \cite[Theorem 3.2.5]{CM-Book}, a property that quasiconvex functions do not have.
%\medskip
%Quasiconvexity is the mathematical formulation of the assumption {\it tendency to diversification} on the consumers in consumer's preference theory (see \cite{D-1959}). 
Furthermore, the following geometric characterization is well-known:
 \begin{align*}
  h ~ \mathrm{is ~ convex} & \Longleftrightarrow \, \mathrm{epi}\,h ~
  \mathrm{is ~ a ~ convex ~ set.}\\ 
  h ~ \mathrm{is ~ quasiconvex} & \Longleftrightarrow \, S_{\delta} (h) ~ \mathrm{is ~ a ~ convex ~ set ~ for ~ all ~ } \delta \in \mathbb{R}.
 \end{align*}

The following definition (see \cite{KL,NTV}) includes (strongly) convex, (strongly) quasiconvex, (strongly) star-convex and (strongly) quasar-convex functions when the set of minimizers is nonempty.

\begin{definition}\label{def:starqc}
 Let $K \subset \mathbb{R}^n$ be a nonempty convex set ($K = \mathbb{R}^n$ unless a proper subset is specified), $h: K \rightarrow \mathbb{R}$ be a function, and $\overline{x} \in {\rm argmin}_{ K}\,h$. Then $h$ is (strongly) star quasiconvex on $K$ with respect to $\overline{x}$ if there exists $\gamma \geq 0$ such that for every $y \in K$ and every $\lambda \in [0, 1]$, we have
 \begin{align}\label{ss:qcx}
  h(\lambda \overline{x} + (1-\lambda)y) \leq h(y) - \lambda (1 - \lambda) \frac{\gamma}{2} \lVert y - \overline{x} \rVert^{2}.
 \end{align}
 It is said that $h$ is strongly star quasiconvex on $K$ when $\gamma > 0$ and it is star quasiconvex on $K$ when $\gamma = 0$, in both cases, with respect to a specific minimizer $\overline{x}$. When $K = \mathbb{R}^n$, we omit reference to $K$ and simply say that $h$ is (strongly) star quasiconvex with respect to $\overline{x}$.
\end{definition}

We note immediately that when $\gamma > 0$, strongly star quasiconvex functions on $K$ have a unique minimizer on $K$. As a consequence, in this case, we simply say that $h$ is strongly star quasiconvex on $K$ with modulus $\gamma > 0$, without re\-fe\-ren\-ce to a point.

Clearly, every (strongly) quasiconvex function is star quasiconvex (on its domain), while the reverse implication does not hold true (see \cite{KL}).

Let us summarize the most interesting/important properties related to star quasiconvexity:

\begin{lemma}\label{already:know}
\begin{itemize} 
 \item[$(a)$] {\bf Star-shaped Sublevel Sets} {\rm (\cite[Theorem 9]{KL})}: Let $h: \mathbb{R}^{n} \rightarrow \overline{ \mathbb{R}}$ be a proper function. Then, $h$ is star quasiconvex with respect to $\overline{x} \in {\rm argmin}_{ \mathbb{R}^{n}}\,h$ if and only if $S_{\delta} (h)$ is star-shaped at $\overline{x}$ for all $\delta \in \mathbb{R}$.

 \item[$(b)$] {\bf Nondecreasing Along Rays} {\rm (\cite[Corollary 15]{KL})}: A function $h$ is star quasiconvex ($\gamma=0$) with respect to $\overline{x} \in {\rm argmin}_{\mathbb{R}^{n}}\,h$ if and only if $h$ is nondecreasing along rays through $\overline{x}$ on $\mathbb{R}^{n}$.

 \item[$(c)$] {\bf Separable Characterization}  {\rm (\cite[Theorem 3]{Lara-10})}: Let $X \subset \mathbb{R}^{n}$ be a convex set such that $X := \prod^{m}_{i=1} X_{i}$, where each $X_{i} \subset \mathbb{R}^{n_{i}}$ is a convex set, $n_{i} \geq 1$ and $\sum^{m}_{i=1} n_{i} = n$. Let $\phi_{i}: X_{i} \rightarrow \mathbb{R}$ be a function ($i = 1,\ldots,m$) and $\phi: X \rightarrow \mathbb{R}$ be a function such that
\begin{equation}\label{separable}
 \phi(x_{1}, \ldots x_{m}) := \sum \limits_{i=1}^{m} \phi_{i}(x_{i}). 
\end{equation}
 Then, $\phi$ is (strongly) star quasiconvex with modulus $\gamma \geq 0$ with respect to $\overline{x} := (\overline{x}_{1}, \ldots, \overline{x}_{m}) \in {\rm argmin}_{X}\phi$ if and only if every $\phi_{i}$ is (strongly) star quasiconvex with modulus $\gamma \geq 0$ with respect to $\overline{x}_{i} \in {\rm argmin}_{X_{i}} \phi_{i}$.
 
 \item[$(d)$] {\bf First-Order Characterization}  {\rm (\cite[Theorem 3.3]{NTV})}:
 Let $h: \mathbb{R}^{n} \rightarrow \mathbb{R}$ be a di\-ffe\-ren\-tia\-ble function and $\overline{x} \in {\rm argmin}_{\mathbb{R}^{n}}\,h$. Then, $h$ is (strongly) star qua\-si\-con\-vex with modulus $\gamma \geq 0$ with respect to $\overline{x}$ if and only if for every $y \in \mathbb{R}^{n}$, we have
 \begin{equation}\label{char:grad}
  \langle \nabla h(y), y - \overline{x} \rangle \geq \frac{\gamma}{2} \lVert y - \overline{x} \rVert^{2}.
 \end{equation}
\end{itemize} 
\end{lemma}

\begin{remark}\label{rem:important}
\begin{itemize}
 \item[$(i)$] {\rm (see \cite[Remark 2$(iii)$]{Lara-10}): Observe that when $\gamma > 0$, part $(d)$ shows that differentiable strongly star quasiconvex functions coincides with functions satisfying the {\it Restricted Secant Inequality} (RSI) property (see \cite{Yi} among others). When $\gamma = 0$, differentiable star quasiconvex functions include {\it Variationally Coherent} (VC) functions, which were developed in \cite{ZBoyd} (or differentiable functions which satisfy \cite[assumption $(1.2)$]{SS}). This inclusion may be strict as \cite[Remark 28$(ii)$]{KL} shows.

 Both classes (RSI) and (VC) have been studied deeply in continuous optimization and machine learning theory in virtue of its properties for gradient-type methods. }

 \item[$(ii)$] {\rm  (see \cite[Corollary 4]{Lara-10}): As a consequence of the characterization in part $(c)$, if each $\phi_{i}$ is quasiconvex and $\overline{x}_{i} \in {\rm argmin}_{X_{i}} \phi_{i}$, then $\phi = \sum_{i} \phi_{i}$ is star quasiconvex with respect to $\overline{x} = (\overline{x}_{1}, \ldots, \overline{x}_{m}) \in \operatorname{argmin}_{X} \phi$, thereby resolving the converse implication of the Debreu--Koopmans theorem \cite{DK} (see also \cite{CL}). This theorem is a useful tool in economic modeling and microeconomics, particularly in production theory and consumer preference theory.}

 \item[$(iii)$] {\rm Note that functions with star-shaped sublevel sets have interesting applications in consumer's preference theory (see \cite{FOS,Mas} and references therein). } 
\end{itemize}
\end{remark}

Let $K \subset \mathbb{R}^n$ be a nonempty set and $T: K \rightrightarrows \mathbb{R}^{n}$ be a set-valued operator. We write $\zer(T) := \{x \in K: \, 0\in T(x)\}$ and $\ran(T) := \{u\in \mathbb{R}^n: \, \exists\, x\in K, \, u \in T(x)\}$. Furthermore, $T$ is said to be:
\begin{itemize}
 \item[$(i)$] (Strongly) monotone with modulus $\gamma \geq 0$ on $K$ if
 \begin{equation}\label{def:mon:set}
  \langle v - w, x - y \rangle \geq \frac{\gamma}{2}\|y - x\|^{2}, ~ \forall ~ x, y \in K, ~ \forall ~ v \in T(x), ~ \forall ~ w \in T(y);
 \end{equation}

 \item[$(ii)$] (Strongly) pseudomonotone with modulus $\gamma \geq 0$ on $K$ if
 \begin{equation}\label{def:pseudo:set}
  \langle u, y - x \rangle \geq 0 ~ \Longrightarrow ~ \langle v, y - x \rangle \geq \frac{\gamma}{2}\|y - x\|^{2}, ~ \forall ~ x, y \in K, ~ \forall ~ u \in T(x), ~ \forall ~ v \in T(y);
 \end{equation}

 \item[$(iii)$] (Strongly) quasimonotone with modulus $\gamma \geq 0$ on $K$ if
 \begin{equation}\label{def:quasi:set}
  \langle u, y - x \rangle > 0 ~ \Longrightarrow ~ \langle v, y - x \rangle \geq \frac{\gamma}{2}\|y - x\|^{2}, ~ \forall ~ x, y \in K, ~ \forall ~ u \in T(x), ~ \forall ~ v \in T(y).
 \end{equation}
\end{itemize} 
Clearly, every (strongly) monotone operator is (strongly) pseudomonotone, and every (strongly) pseudomonotone operator is (strongly) quasimonotone. The reverse statement does not hold true in general (see, for instance, \cite{ADSZ,HKS}). Furthermore, note that the pseudomonotonicity notion in \eqref{def:pseudo:set} is the generalized monotonicity notion of Karamardian \cite{Kara}, rather than the one due to Brezis \cite{Brezis1968}.

\medskip

The previous notions are related to the generalized convexity theory in virtue of the following results. Let $h: \mathbb{R}^{n} \to \mathbb{R}$ be a differentiable function. Then the following assertions hold.
\begin{itemize}
 \item[$(i)$] $h$ is (strongly) convex with modulus $\gamma \geq 0$ if and only if $\nabla h$ is (strongly) monotone with modulus $2\gamma \geq 0$.

 \item[$(ii)$] $h$ is pseudoconvex if and only if $\nabla h$ is pseudomonotone.

 \item[$(iii)$] $h$ is quasiconvex if and only if $\nabla h$ is quasimonotone.
\end{itemize}

Now, let us recall some continuity notions for set-valued operators.

\begin{definition}
 Let $K \subset \mathbb{R}^n$ be a nonempty set and  $T: K \rightrightarrows \mathbb{R}^n$ be a set-valued operator. A function $t: K \to \mathbb{R}^n$ is called a {\it selection} of $T$ when
 \[
  t(x) \in T(x), ~ \forall ~ x \in K.
 \]
\end{definition}

\begin{definition}\label{def:hemi}
 Let $K \subset \mathbb{R}^n$ be an open set and $T: K \rightrightarrows \mathbb{R}^n$ be a set-valued operator. Then $T$ is said to be:
 \begin{itemize}
  \item[$(i)$] hemicontinuous at $x \in K$ when for every $w \in T(x)$ and every $d \in \mathbb{R}^n$, there exist $\varepsilon > 0$ and a map $t: \, ]0, \varepsilon[ \, \to \mathbb{R}^n$ such that $x + s d \in K$ and $t(s) \in T(x + s d)$ for all $s \in \, ]0, \varepsilon[$, and
  \[
   \lim_{s \to 0^+} t(s) = w;
  \]
  $T$ is hemicontinuous on $K$ when it is hemicontinuous at every $x \in K$;

  \item[$(ii)$] weakly hemicontinuous at $x_0 \in K$ when for every $w \in T(x_0)$ and every $d \in \mathbb{R}^n$, there exist $\varepsilon > 0$ and a map $t: \, ]0, \varepsilon[ \, \to \mathbb{R}^n$ such that $x_0 + s d \in K$ and $t(s) \in T(x_0 + s d)$ for all $s \in \, ]0, \varepsilon[$, and
  \[
   \lim_{s \to 0^+} \langle t(s), w \rangle = \langle w, w \rangle = \|w\|^2.
  \]
 \end{itemize}
\end{definition}

\noindent Since $K$ is open, for every $x \in K$ and $d \in \mathbb{R}^n$ there exists $\varepsilon > 0$ such that $x + s d \in K$ for all $s \in \, ]0, \varepsilon[$. Clearly, hemicontinuity at $x_0$ implies weak hemicontinuity at $x_0$ (take the inner product with $w$), while the reverse statement does not hold in general (see \cite{BI-book, HKS}). When $T$ is single-valued, hemicontinuity at $x$ reduces to the classical condition $\lim_{s \to 0^+} T(x + s d) = T(x)$ for every $d \in \mathbb{R}^n$.

%Note that hemicontinuity implies weak hemicontinuity, while the reverse statement does not hold in general (see \cite{BI-book, HKS}).

For a further study on monotone operators, generalized monotonicity and generalized convexity, we refer to \cite{ADSZ,BI-book,CM-Book,HKS,Lara-9,Lara-10} and references therein.

\section{Star Monotone Operators}\label{sec:03}

As explained in the introduction, we are interested in the study of the following generalized monotonicity notion, called {\it star monotonicity}. It is important to mention that this notion is not new, since it has appeared as an assumption of sufficient condition for ensuring the convergence of algorithms for variational inequalities (see, for instance, \cite[equation (1.2)]{SS} and also \cite{CGL2023,KLL2022,Mal}).

\begin{definition}\label{star:mon}
 Let $K \subset \mathbb{R}^{n}$ be a nonempty set and $T: K \rightrightarrows \mathbb{R}^{n}$ be a set-valued operator. Then $T$ is said to be {\it (strongly) star monotone} with modulus $\gamma \geq 0$ at $x_{0} \in K$ when
 \begin{equation}\label{eq:star:mon}
   \left \langle v, y - x_{0} \right \rangle \geq \frac{\gamma}{2}\left \Vert y - x_{0} \right \Vert ^{2}, ~ \forall ~ v \in T(y), ~ \forall ~ y \in K.
 \end{equation}
 It is called strongly star monotone when $\gamma > 0$, and star monotone when $\gamma = 0$. We refer to both notions as ``(strongly) star monotone with modulus $\gamma \geq 0$".
\end{definition}

Observe that every (strongly) monotone operator with $\zer (T) \neq \emptyset$ is (strong\-ly) star monotone at any point $x_{0} \in \zer (T)$. The converse statement does not hold as the following remark shows.

\begin{remark}\label{no:rel:mon}
 Let $f: \mathbb{R}^{2} \rightarrow \mathbb{R}$ be defined by $f(x_{1}, x_{2}) = f_{1} (x_{1}) + f_{2} (x_{2})$, with $f_{i} (x_{i}) = \frac{x^{2}_{i}}{1 + x^{2}_{i}}$ for $i=1,2$. Note that $f$ is differentiable and nonconvex, thus $\nabla f$ is not monotone. Since every $f_{i}$ is quasiconvex without being convex, $f=f_{1} + f_{2}$ is neither quasiconvex nor pseudoconvex by the Debreu--Koopmans theorem (see \cite[Theorems 2 and 10]{DK}), thus $\nabla f$ is neither quasimonotone nor pseudomonotone.
 
 On the other hand, since $0 \in {\rm argmin}_{ \mathbb{R}}\,f_{i}$ and $f_{i}$ is quasiconvex for $i=1,2$, $f=f_{1}+f_{2}$ is star quasiconvex with respect to $(0, 0) \in {\rm argmin}_{\mathbb{R}^{2}}\,f$ by \cite[Corollary 4]{Lara-10} (see Remark \ref{rem:important}$(ii)$), and it follows from \cite[Theorem 3.3]{NTV} (see Lemma \ref{already:know}$(d)$) that $\nabla f$ is star monotone at $(0, 0)$.
\end{remark}

\begin{remark}\label{no:rel:mon:robust}
 As a second illustration in a different applied context, consider the {\it Cauchy} (or Lorentzian) loss function, widely used in robust regression, computer vision, and robust M-estimation: $g_i(x_i) := \log(1+x_i^2)$ for $i=1,2$, and $g(x_1,x_2) := g_1(x_1) + g_2(x_2)$. Each $g_i$ is differentiable, even, with unique minimizer at $0$ (since $g_i'(x_i) = \frac{2x_i}{1+x_i^2}$ is negative for $x_i<0$ and positive for $x_i>0$), hence quasiconvex on $\mathbb{R}$; it is not convex, since $g_i''(x_i) = \frac{2(1-x_i^2)}{(1+x_i^2)^2} < 0$ for $|x_i|>1$.

 Since $\log$ is increasing, $g(x_1,x_2) = \log\big[ (1+x_1^2)(1+x_2^2) \big]$ is quasiconvex if and only if $h(x_1,x_2) := (1+x_1^2)(1+x_2^2)$ is quasiconvex. Taking $P_1=(3,0)$ and $P_2=(0,3)$, we have $h(P_1)=h(P_2)=10$, while $h(1.5,1.5) = (1+2.25)^2 = 10.5625 > 10$, directly violating the quasiconvexity inequality $h\big( \tfrac{1}{2} P_1 + \tfrac12P_2\big) \leq \max\{h(P_1), h(P_2)\}$. Hence $g$ is not quasiconvex, and since differentiable pseudoconvex functions are always quasiconvex, $g$ is not pseudoconvex either; consequently $\nabla g$ is neither quasimonotone nor pseudomonotone.

 On the other hand, $(0,0) \in {\rm argmin}_{ \mathbb{R}^2} g$, and since each $g_i$ is quasiconvex with minimizer $0$, $g$ is star quasiconvex with respect to $(0,0)$ by Lemma \ref{already:know}$(c)$ (\cite[Corollary 4]{Lara-10}), so $\nabla g$ is star monotone at $(0,0)$ by Lemma \ref{already:know}$(d)$ (\cite[Theorem 3.3]{NTV}), despite failing every classical generalized monotonicity notion, exactly as in Remark \ref{no:rel:mon}, but now in a setting drawn from robust statistics rather than an abstract construction.
\end{remark}

\subsection{Basic Properties}

As a first result, we present the usual calculus rules. The easy proof is omitted.

\begin{proposition}\label{prop:01:set}
 Let $K \subset \mathbb{R}^{n}$ be a nonempty set and $T: K \rightrightarrows \mathbb{R}^n$ be a (strongly) star monotone at $x_0 \in K$ with modulus $\gamma \geq 0$. Then the following assertions hold.
 \begin{itemize}
  \item[$(a)$] For every $\alpha > 0$, the operator $\alpha T$ is (strongly) star monotone at $x_0$ with modulus $\alpha \gamma \geq 0$.

  \item[$(b)$] For every $x^\ast \in K$, if $\widetilde{T}: (K - x^\ast) \rightrightarrows \mathbb{R}^n$ is defined by $\widetilde{T}(z) := T(z+x^\ast)$ (so that $z + x^\ast \in K$), then $T$ is (strongly) star monotone at $x^\ast$ with modulus $\gamma \geq 0$ if and only if $\widetilde{T}$ is (strongly) star monotone at $0$ with the same modulus $\gamma$. In particular, taking $x^\ast = x_0$ recovers a translation of the hypothesis above. % parts $(b)$ and $(c)$ hold for an arbitrary anchor point and do not use the specific hypothesis that $T$ is star monotone at $x_0$.)

  \item[$(c)$] For every $x^\ast \in K$, if $T$ has any of the following properties: lower semicontinuity, nonempty values, closed values, convex values, or closed graph, then $\widetilde{T}$ (as defined in $(b)$) has the same property.

  \item[$(d)$] For every $y \in K$ and every $v \in T(y)$, we have
  \[
   \|v\| \geq \frac{\gamma}{2}\|y - x_0\|.
  \]
 \end{itemize}
\end{proposition}

\begin{proof}
 $(a)$: Directly from the definition, $(b)$ follows from $\langle u, z\rangle = \langle u, (z+x^\ast)-x^\ast \rangle$ with $u \in T(z+x^\ast)$, and $(c)$ is immediate from the definition of $\widetilde{T}$. 
 
 $(d)$: From (strong) star monotonicity with modulus $\gamma \geq 0$ and Cauchy-Schwarz inequality, we have
 \[
  \|v\|\|y - x_0\| \geq \langle v, y - x_0 \rangle \geq \frac{\gamma}{2}\|y - x_0\|^2.
 \]
 If $y \neq x_0$, then we divide by $\|y - x_0\|$; if $y = x_0$, then the required inequalities are trivial.
\end{proof}

As a consequence of Proposition \ref{prop:01:set}$(d)$, we have the following:

\begin{corollary} \label{cor:norm:consequences}
 Let $K \subset \mathbb{R}^n$ be a nonempty set and $T: K \rightrightarrows \mathbb{R}^n$ be strongly star monotone with modulus $\gamma > 0$ at $x_0 \in K$. Then the following assertions hold (with the convention $\inf \emptyset = +\infty$).
 \begin{itemize}
  \item[$(a)$] If $K$ is unbounded, then $T$ is coercive, i.e.,
  \[
   \liminf_{y \in K,\, \|y\| \to+ \infty} \, \inf_{v \in T(y)} \frac{\|v\|}{\|y\|} \geq \frac{\gamma}{2} > 0.
  \]

  \item[$(b)$] The radial component of $T$ grows at least linearly, i.e.,
  \[
   \left\langle v, \frac{y - x_0}{\|y - x_0\|} \right\rangle \geq \frac{\gamma}{2}\|y - x_0\|, ~ \forall ~ y \in K \setminus \{x_0\}, ~ \forall ~ v \in T(y).
  \]
  In particular, if $x_0$ is an accumulation point of $K$, then
  \[
   \liminf_{y \in K \setminus \{x_0\},\, y \to x_0} \, \inf_{v \in T(y)} \frac{\langle v, y - x_0 \rangle}{\|y - x_0\|^2} \geq \frac{\gamma}{2}.
  \]

  \item[$(c)$] For every $c \geq 0$,
  \[
   L_c := \{y\in K : \exists\, v\in T(y),\ \|v\|\leq c\} \subset \overline{\mathbb{B}} (x_0, 2c/\gamma).
  \]
  In particular, $\zer(T) \subset L_0 \subset \{x_0\}$.
 \end{itemize}
\end{corollary}

\begin{proof}
 $(a)$: By Proposition \ref{prop:01:set}$(d)$, every $v\in T(y)$ satisfies $\|v\| \geq \frac{\gamma}{2}\|y-x_0\|$, and this bound holds trivially when $T(y) = \emptyset$ by the convention $\inf \emptyset = +\infty$. Dividing by $\|y\|>0$, we obtain
 \[
  \inf_{v\in T(y)} \frac{\|v\|}{\|y\|} \geq \frac{\gamma}{2} \frac{\|y-x_0\|}{\|y\|} \geq \frac{\gamma}{2} \left( 1 - \frac{\|x_0\|}{\|y\|} \right), ~ \forall ~ y \in K \setminus \{0\}.
 \]
 Since $K$ is unbounded, the result follows by taking the limit inferior as $\|y\| \to +\infty$.

 $(b)$: For every $y \in K \setminus \{x_0\}$ and every $v \in T(y)$, strong star monotonicity gives $\langle v, y - x_0 \rangle \geq \frac{\gamma}{2}\|y - x_0\|^2$. Dividing by $\|y - x_0\| > 0$ yields the first inequality, while dividing by $\|y - x_0\|^2 > 0$, taking the infimum over $v \in T(y)$, and then the limit inferior as $y \to x_0$ yields the second one.

 $(c)$: If $y \in L_c$, there exists $v\in T(y)$ with $\|v\|\leq c$. By Proposition \ref{prop:01:set}$(d)$,
 \[
  c \geq \|v\| \geq \frac{\gamma}{2}\|y-x_0\|,
 \]
 hence $\|y-x_0\| \leq 2c/\gamma$, i.e., $y \in \overline{\mathbb{B}} (x_0, 2c/\gamma)$. For $c = 0$, this gives $L_0 \subset \{x_0\}$, and clearly $\zer(T) \subset L_0$.
\end{proof}

\begin{proposition}\label{prop:02:set}
 Let $K \subset \mathbb{R}^{n}$ be a nonempty set  and $T_1, \dots, T_m: K \rightrightarrows \mathbb{R}^n$ be (strongly) star monotone at $x_0 \in K$ with modulus $\gamma_1, \dots, \gamma_m \geq 0$. Then for any $\alpha_1, \dots, \alpha_m \geq 0$, the operator $A: K \rightrightarrows \mathbb{R}^n$ defined by
 \begin{equation}\label{sum:weights}
  A(y) := \sum_{i=1}^m \alpha_i T_i(y) = \left\{ \sum_{i=1}^m \alpha_i v_i : v_i \in T_i(y) \right\}
 \end{equation}
 is (strongly) star monotone at $x_0$ with modulus $\sum_{i=1}^m \alpha_i \gamma_i$.
\end{proposition}

\begin{proof}
 For any $y \in K$ and any $v \in A(y)$, there exist $v_i \in T_i(y)$ such that $v = \sum_{i=1}^m \alpha_i v_i$. Then,
 \begin{align*}
  \langle v, y - x_0 \rangle = \sum_{i=1}^m \alpha_i \langle v_i, y - x_0 \rangle & \geq \sum_{i=1}^m \alpha_i \frac{\gamma_i}{2} \|y - x_0\|^2 = \frac{1}{2} \left( \sum_{i=1}^m \alpha_i \gamma_i \right) \|y - x_0\|^2,
 \end{align*}
 and the result follows.
\end{proof}

\begin{proposition}\label{prop:03:set}
 Let $K \subset \mathbb{R}^{n}$ be a nonempty set, $T: K \rightrightarrows \mathbb{R}^n$ be (strongly) star monotone at $x_0 \in K$ with modulus $\gamma \geq 0$, $A: \mathbb{R}^m \to \mathbb{R}^n$ be a linear operator with adjoint $A^*$, and $b \in \mathbb{R}^n$. Set
 \[
  K_A := \{ y \in \mathbb{R}^m : \, Ay + b \in K \},
 \]
 and define $S: K_A \rightrightarrows \mathbb{R}^m$ by
 \[
  S(y) := A^* T(Ay + b) = \{ A^* v : \, v \in T(Ay + b) \}.
 \]
 If $y_0 \in \mathbb{R}^m$ satisfies $Ay_0 + b = x_0$ (so that $y_0 \in K_A$), then $S$ is (strongly) star monotone at $y_0$ with modulus $\gamma \lambda_{\min}(A^*A) \geq 0$, where $\lambda_{\min}(A^*A) \geq 0$ denotes the smallest eigenvalue of $A^*A$. In particular, if $\gamma > 0$ and $A$ is injective, then $S$ is strongly star monotone at $y_0$ with modulus $\gamma \lambda_{\min}(A^*A) > 0$.
\end{proposition}

\begin{proof}
 Let $y \in K_A$ and $w \in S(y)$. Then there exists $v \in T(Ay + b)$ such that $w = A^* v$. Since $Ay + b \in K$ and $(Ay + b) - x_0 = A(y - y_0)$, the (strong) star monotonicity of $T$ at $x_0$ yields
 \begin{align*}
  \langle w, y - y_0 \rangle
  & = \langle A^* v, y - y_0 \rangle = \langle v, A(y - y_0) \rangle = \langle v, (Ay + b) - x_0 \rangle \\
  & \geq \frac{\gamma}{2}\|A(y - y_0)\|^2 = \frac{\gamma}{2} \langle A^*A(y - y_0), y - y_0 \rangle \\
  & \geq \frac{\gamma}{2} \lambda_{\min}(A^*A) \|y - y_0\|^2,
 \end{align*}
 where the last inequality follows from the Rayleigh quotient bound for the symmetric positive semidefinite operator $A^*A$. Hence, $S$ is (strongly) star monotone at $y_0$ with modulus $\gamma \lambda_{\min}(A^*A)$.

 Finally, $\lambda_{\min}(A^*A) = \min_{\|z\| = 1} \|Az\|^2$, which is positive if and only if $\ker A = \{0\}$, i.e., if and only if $A$ is injective. This proves the last assertion.
\end{proof}

In the case of strongly star monotone operators ($\gamma > 0$), we can obtain further interesting properties.

\begin{proposition} \label{prop:04}
 Let $K \subset \mathbb{R}^{n}$ be a nonempty set  and $T: K \rightrightarrows \mathbb{R}^{n}$ be a (strongly) star monotone with modulus $\gamma > 0$ at $x_0 \in K$, $\varepsilon > 0$ and $S: K \to \mathbb{R}^n$ be such that $\|S(y)\| \leq \varepsilon$ for all $y \in K$. If $\varepsilon < \gamma \frac{\delta}{2}$ for some $\delta > 0$, then for every $y \in K \setminus \mathbb{B} (x_0, \delta)$ and every $v \in T(y)$, we have
 \[
  v + S(y) \neq 0.
 \]
 In other words, $0 \notin (T + S)(y)$ for all $y \in K \setminus \mathbb{B} (x_0, \delta)$.
\end{proposition}

\begin{proof}
 For any $y \in K$ with $\|y - x_0\| \geq \delta$ and any $v \in T(y)$, we have
 \begin{align*}
  \langle v + S(y), y - x_0 \rangle & \geq \frac{\gamma}{2}\|y - x_0\|^2 - \|S(y)\|\|y - x_0\| \\
  & \geq \left(\frac{\gamma}{2}\|y - x_0\| - \varepsilon \right) \|y - x_0\| \geq \left( \frac{\gamma\delta}{2} - \varepsilon \right) \|y - x_0\| > 0.
 \end{align*}
 Hence, $\langle v + S(y), y - x_0 \rangle > 0$, which implies $v + S(y) \neq 0$. Therefore, $T + S$ has no zeros outside the ball $\mathbb{B} (x_0, \delta)$.
\end{proof}

We finish this subsection with the following two remarks, which show properties that (strongly) star monotone operators do not have.

\begin{remark} \label{rem:no:convex}
 Note that $T(y)$ need not be convex, closed, or star-shaped at any point, even when $T$ is strongly star monotone with modulus $\gamma > 0$ at $x_{0}$. Indeed, let $K = \mathbb{R}$, $x_0 = 0$, $\gamma = 2$, and define
  \[
  T(y) := 
  \begin{cases} 
   \{2y + 1\} \cup \{y + 1/n : n \in \mathbb{N}\} & {\rm if } ~~ y > 0, \\ 
   \{y\} & {\rm if } ~~ y \leq 0. 
  \end{cases}
 \]
 
 For $y > 0$, every $v \in T(y)$ satisfies that $v \geq y > 0$ and $\langle v, y\rangle = vy \geq y^2 = \frac{\gamma}{2}y^2$ for all $v \in T(y)$; for $y \leq 0$, every $v \in T(y)$ satisfies that $\langle v,y\rangle = y^2 = \frac{\gamma}{2}y^2$. Therefore, $T$ is strongly star monotone at $0$ with modulus $\gamma=2$.

 Note that for $y>0$, $T(y)$ is clearly nonconvex and non star-shaped at any point (since it has isolated points).
 
 Furthermore, for $y>0$, $T(y)$ is not closed. Indeed, take $\{v_{n} = y + \frac{1}{n}\}_{n} \subset T(y)$, then $y+1/n \to y \not\in T(y)$, hence $T(y)$ is not closed.
\end{remark}

\begin{remark} \label{rem:unbounded}
 $T(y)$ is not necessarily bounded, even under strong star monotonicity with $\gamma > 0$. Indeed, let $K = [0, + \infty[ \, \subset \mathbb{R}$, $x_0 = 0$, and fix any $\gamma \in \, ]0,2]$. Define
 \[
  T(y) := 
  \begin{cases} 
   [\,y, +\infty[ & {\rm if } ~~ y > 0, \\ 
   \{0\} & {\rm if } ~~ y = 0. 
  \end{cases}
 \]
 Clearly, $T(y)$ is unbounded for all $y>0$.
 
 Furthermore, $T$ is strongly star monotone at $x_0=0$ with modulus $\gamma \in \, ]0, 2]$. Indeed, for $y=0$, it is trivial. For $y>0$ and any $v \in T(y) = [y, + \infty[$, we have $v \geq y > 0$, so $\langle v, y - x_0 \rangle = vy \geq  y^2 \geq \frac{\gamma}{2} y^2 = \frac{\gamma}{2}\|y-x_0\|^2$, using $\gamma \leq 2$ in the last inequality. Therefore, $T$ is strongly star monotone at $0$ with modulus $\gamma \in \, ]0, 2]$. 
\end{remark}

\subsection{Existence of Zeroes}

Star monotone operators are closely related to the existence of zeroes. Indeed, as we show in the next result, if a weakly hemicontinuous operator $T$ satisfies relation \eqref{eq:star:mon} at a point $x_{0}$, then $x_{0}$ should be a zero of $T$. We emphasize that in the following result, no monotonicity assumption is needed.

\begin{theorem}\label{theo:zero}
 Let $K \subset \mathbb{R}^{n}$ be an open set and $T: K \rightrightarrows \mathbb{R}^{n}$ be a weakly hemicontinuous set-valued operator at $x_0 \in K$ with $T(x_{0}) \neq \emptyset$. Then the following assertions hold.
 \begin{itemize}
  \item[$(a)$] If $T$ is (strongly) star monotone at $x_0$ with modulus $\gamma \geq 0$, then $0 \in T(x_0)$.

  \item[$(b)$] If $T$ is strongly star monotone with modulus $\gamma > 0$ at $x_0 \in K$, then $x_0$ is the unique point in $K$ satisfying \eqref{eq:star:mon}. As a consequence, $\zer (T) = \{x_0\}$, that is, $x_0$ is the unique zero of $T$.
 \end{itemize} 
\end{theorem}

\begin{proof}
 Suppose for the contrary that $0 \notin T(x_0)$, and choose $w \in T(x_0)$, so $w \neq 0$. By weak hemicontinuity at $x_0$ with direction $d = -w$, there exist $\varepsilon > 0$ and $t: \, ]0, \varepsilon[ \, \to \mathbb{R}^n$ with $y_s := x_0 - s w \in K$, $t(s) \in T(y_s)$ and $\lim_{s \to 0^+} \langle t(s), w \rangle = \|w\|^2$. By (strong) star monotonicity, $\langle t(s), y_s - x_0 \rangle \geq \frac{\gamma}{2}\|y_s - x_0\|^2$, and since $y_s - x_0 = -s w$,
 \[
  - \langle t(s), w \rangle \geq \frac{\gamma}{2} s \|w\|^2, ~ \forall ~ s \in \, ]0, \varepsilon[.
 \]
 Taking the limit as $s \to 0^+$ yields $-\|w\|^2 \geq 0$, a contradiction. Therefore, $0 \in T(x_0)$.

 $(b)$: Since $T$ is weakly hemicontinuous at $x_0$ and strongly star monotone at $x_0$ with modulus $\gamma > 0$, part $(a)$ yields $0 \in T(x_0)$. 
 
 Suppose that $x_1 \in K$ also satisfies
 \begin{equation}\label{ineq:01:set}
  \langle v, y - x_1 \rangle \geq \frac{\gamma}{2}\|y - x_1\|^2, ~ \forall ~ v \in T(y), ~ \forall ~ y \in K.
 \end{equation}
 Taking $y = x_0$ and $v = 0 \in T(x_0)$ in \eqref{ineq:01:set}, we obtain
 \[
  0 = \langle 0, x_0 - x_1 \rangle \geq \frac{\gamma}{2}\|x_0 - x_1\|^2.
 \]
 Since $\gamma > 0$, we have $\|x_0 - x_1\|^2 \leq 0$, and hence $x_1 = x_0$.

 Suppose for contradiction that there exists $x_1 \in \zer (T)$ with $x_1 \neq x_0$. Then $0 \in T(x_1)$. Since $T$ is strongly star monotone at $x_0$, taking $y = x_1$ and $v = 0 \in T(x_1)$ in \eqref{eq:star:mon} gives $0 = \langle 0, x_1 - x_0 \rangle \geq \frac{\gamma}{2}\|x_1 - x_0\|^2$. Since $\gamma > 0$ and $x_1 \neq x_0$, we have $\|x_1 - x_0\|^2 > 0$, a contradiction. Therefore, $\zer (T)$ is a singleton.
\end{proof}

Let us now present the following result, which shows that strong star monotonicity ($\gamma > 0$) is global metric subregularity.

\begin{corollary}[Global Metric Subregularity] \label{cor:metric:subreg}
 Let $K \subset \mathbb{R}^n$ be an open set and $T: K \rightrightarrows \mathbb{R}^n$ be weakly hemicontinuous at $x_0 \in K$, strongly star monotone at $x_0$ with modulus $\gamma > 0$, and $T(x_{0}) \neq \emptyset$. Then $T$ is metrically subregular at $x_0$ for $0$ with modulus $\kappa = 2/\gamma$, globally on $K$, that is,
 \begin{equation}\label{global:subregu}
  \|y - x_0\| \leq \frac{2}{\gamma} {\rm dist} (0, T(y)), ~ \forall ~ y \in K,
 \end{equation}
 where ${\rm dist} (0, T(y)) := \inf_{v \in T(y)} \|v\|$.
\end{corollary}

\begin{proof}
 By Theorem \ref{theo:zero}$(b)$, $\zer(T) = \{x_0\}$, thus ${\rm dist} (y, T^{-1}(0)) = \|y-x_0\|$ for every $y\in K$. By Proposition \ref{prop:01:set}$(d)$, every $v\in T(y)$ satisfies $\|v\| \geq \frac{\gamma}{2} \|y-x_0\|$. Taking the infimum over $v\in T(y)$ gives ${\rm dist} (0, T(y)) \geq \frac{\gamma}{2}\|y-x_0\|$, i.e., $\|y-x_0\| \leq \frac{2}{\gamma} {\rm dist} (0, T(y))$.
\end{proof}

\begin{remark}
 Unlike the classical notion, which is a local property, Corollary \ref{cor:metric:subreg} gives subregularity on all of $K$, with an explicit modulus $\kappa=2/\gamma$ determined entirely by the star monotonicity modulus. This places strongly star monotone operators inside the broader metric regularity framework (see, e.g., \cite{DR}) that underlies most modern linear-convergence proofs for first-order methods. 
\end{remark}

Another consequence is the following Tikhonov-type regularization result.

\begin{corollary} \label{cor:tikhonov}
 Let $K \subset \mathbb{R}^n$ be open, $x_0 \in K$, and $T: K \rightrightarrows \mathbb{R}^n$ be star monotone at $x_0$ with modulus $\gamma \geq 0$ with $T(x_{0}) \neq \emptyset$. For every $\varepsilon > 0$, define
 \begin{equation}\label{tikho:reg}
  T_\varepsilon(y) := T(y) + \varepsilon(y - x_0).
 \end{equation}
 Then $T_\varepsilon$ is strongly star monotone at $x_0$ with modulus $\gamma + 2\varepsilon > 0$. Consequently, if $T_\varepsilon$ is weakly hemicontinuous at $x_0$, then $\zer( T_{\varepsilon}) = \{x_0\}$.
\end{corollary}

\begin{proof}
 The operator $S(y) := y - x_0$ satisfies $\langle S(y), y - x_0\rangle = \|y-x_0\|^2 = \frac{2}{2}\|y-x_0\|^2$, i.e., $S$ is strongly star monotone at $x_0$ with modulus $\gamma = 2$. By Proposition \ref{prop:02:set} applied to $T$ and $S$ with weights $\alpha_1 = 1$ and $\alpha_2 = \varepsilon$, the operator $T_\varepsilon = T + \varepsilon S$ is strongly star monotone at $x_0$ with modulus $\gamma + 2\varepsilon$. 
 
 The last claim follows from Theorem \ref{theo:zero}$(b)$.
\end{proof}

%\begin{remark}
% In the set-valued case, the error bound from Proposition \ref{prop:01:set}$(d)$,
%$$\|v\| \geq \frac{\gamma}{2}\|y - x_0\|, ~ \forall ~ v \in T(y),$$
% is the most useful consequence of strong star monotonicity, as it provides a global lower bound on the norm of every selection of $T$ in terms of the distance from $x_0$.
%\end{remark}

\subsection{Relation with Generalized Monotonicity} \label{subsec:3-3}

In this subsection, we analyze the relationship between (strongly) star monotone operators and classical notions in generalized monotonicity theory (monotone, pseudomonotone, and quasimonotone).

% -------------------------------------------------------------------------
\subsubsection{Pseudomonotone Operators}

We begin by establishing that (strong) pseudomonotonicity implies (strong) star monotonicity whenever a zero exists.

\begin{proposition}\label{pseudo:star:set}
Let $K \subset \mathbb{R}^{n}$ be an open set, $T: K \rightrightarrows \mathbb{R}^n$ be a set-valued operator, and $x_0 \in K$. If $0 \in T(x_0)$ and $T$ is (strongly) pseudomonotone with modulus $\gamma \geq 0$ on $K$, then $T$ is (strongly) star monotone at $x_0$ with the same modulus $\gamma \geq 0$.
\end{proposition}

\begin{proof}
Since $0 \in T(x_0)$, for every $y \in K$ and every $u \in T(x_0)$, we have $\langle u, y - x_0 \rangle = 0 \geq 0$. By (strong) pseudomonotonicity of $T$ on $K$ applied at $x = x_0$, every $v \in T(y)$ satisfies $\langle v, y - x_0 \rangle \geq \frac{\gamma}{2}\|y - x_0\|^2$ for all $y \in K$, which proves that $T$ is (strongly) star monotone at $x_0$.
\end{proof}

The converse of Proposition \ref{pseudo:star:set} does not hold in general since the gradient of the function $f$ in Remark \ref{no:rel:mon} is not pseudomonotone.

% -------------------------------------------------------------------------
\subsubsection{Quasimonotone Operators}

Next, we study the relationship between quasimonotone and star monotone operators. We start with a local-to-global result.

\begin{proposition} \label{local:global:set}
 Let $K \subset \mathbb{R}^{n}$ be a convex set and $T: K \rightrightarrows \mathbb{R}^n$ be strongly quasimonotone with modulus $\gamma > 0$ and with nonempty values. If $T$ is strongly star monotone at $x_0 \in K$ on a neighborhood $U \subset K$ of $x_0$, then $T$ is strongly star monotone at $x_0$ on entire $K$.
\end{proposition}

\begin{proof}
 Fix $y \in K \setminus \{x_0\}$, and choose $t_0 \in \, ]0, 1[$ such that $y_t := x_0 + t(y - x_0) \in U$ for all $t \in \, ]0, t_0]$. This is always possible because $K$ is convex and $U$ is a neighborhood of $x_0$.

 For any $t \in \, ]0, t_0]$ and any $u_t \in T(y_t)$, local strong star monotonicity on $U$ yields $\langle u_t, y_t - x_0 \rangle \geq \frac{\gamma}{2}\|y_t - x_0\|^2$. Since $y_t - x_0 = t(y - x_0)$, we get
 \[
  \langle u_t, y - x_0 \rangle \geq \frac{\gamma}{2} t \|y - x_0\|^2 > 0,
 \]
 and consequently $\langle u_t, y - y_t \rangle = (1 - t) \langle u_t, y - x_0 \rangle > 0$.

 By the strong quasimonotonicity of $T$ on $K$, every $v \in T(y)$ satisfies $\langle v, y - y_t \rangle \geq \frac{\gamma}{2}\|y - y_t\|^2$. Taking the limit as $t \to 0^+$ (so $y_t \to x_0$), we obtain
 \[
  \langle v, y - x_0 \rangle \geq \frac{\gamma}{2}\|y - x_0\|^2, ~ \forall ~ v \in T(y).
 \]
 Since $y \in K \setminus \{x_0\}$ was arbitrary, $T$ is strongly star monotone at $x_0$ on entire $K$.
\end{proof}

\begin{remark} \label{not:dropped}
 The strong quasimonotonicity condition ($\gamma > 0$) in Proposition \ref{local:global:set} cannot be dropped, even for single-valued continuous operators (e.g., $T(t) = t$ for $|t| \leq 1$ and $T(t) = 1/t$ for $|t| > 1$).
\end{remark}

\subsubsection{Ray-Wise Characterization}

Finally, motivated by \cite[Theorem 12]{KL}, we characterize star monotonicity via quasimonotonicity along rays emanating from $x_0$.

\begin{proposition} \label{prop:ray:necessity}
 Let $K \subset \mathbb{R}^n$ be a nonempty set, $x_0 \in K$, and $T: K \rightrightarrows \mathbb{R}^n$ be (strongly) star monotone at $x_0$ with modulus $\gamma \geq 0$. Then for every $d \in \mathbb{R}^n \setminus \{0\}$, the restriction of $T$ to the ray $R_d := \{x_0 + sd : s \geq 0\} \cap K$ is (strongly) quasimonotone with modulus $\gamma$.
\end{proposition}

\begin{proof}
 Fix $d \in \mathbb{R}^n \setminus \{0\}$ and pick $x, y \in R_d$. If $x = y$, the statement holds trivially. If $x \neq y$, write $x = x_0 + s_1 d$ and $y = x_0 + s_2 d$ with $s_1, s_2 \geq 0$ and $s_1 \neq s_2$. We consider two cases.

 If $s_1 > s_2$, then $x - x_0 = \frac{s_1}{s_1 - s_2}(x - y)$. Thus, for any $u \in T(x)$, star monotonicity at $x_0$ yields $\langle u, x - x_0 \rangle \geq \frac{\gamma}{2}\|x - x_0\|^2 \geq 0$. Consequently,
 \[
  \langle u, y - x \rangle = -\frac{s_1 - s_2}{s_1} \langle u, x - x_0 \rangle \leq 0,
 \]
 so the premise $\langle u, y - x \rangle > 0$ of (strong) quasimonotonicity is never satisfied, making the implication vacuously true.

 If $s_1 < s_2$, then $y - x = \frac{s_2 - s_1}{s_2}(y - x_0)$. For every $v \in T(y)$, star monotonicity of $T$ at $x_0$ gives $\langle v, y - x_0 \rangle \geq \frac{\gamma}{2}\|y - x_0\|^2$. Multiplying both sides by $\frac{s_2 - s_1}{s_2} \in \, ]0, 1]$ yields
 \begin{align*}
  \langle v, y - x \rangle = \frac{s_2 - s_1}{s_2} \langle v, y - x_0 \rangle & \geq \frac{s_2 - s_1}{s_2} \frac{\gamma}{2} \|y - x_0\|^2 \\
  & = \frac{s_2}{s_2 - s_1} \frac{\gamma}{2} \|y - x\|^2 \geq \frac{\gamma}{2} \|y - x\|^2.
 \end{align*}
 Hence, $T$ restricted to $R_d$ is (strongly) quasimonotone with modulus $\gamma \geq 0$.
\end{proof}

\begin{remark} \label{converse:fails}
 The converse of Proposition \ref{prop:ray:necessity} fails in general (e.g., the constant operator $T(x) \equiv c \neq 0$ is monotone, hence quasimonotone on every ray, but it is not star monotone at any point).
\end{remark}

The converse of Proposition \ref{prop:ray:necessity} holds for $\gamma > 0$ under a radial nondegeneracy condition, which, as we show below, is also necessary.

\begin{theorem}\label{theo:quasi:star:ray}
 Let $K \subset \mathbb{R}^{n}$ be a set star-shaped at $x_0 \in K$, $T: K \rightrightarrows \mathbb{R}^{n}$ be a set-valued operator with nonempty values, and $\gamma > 0$. Assume that $0 \in T(x_0)$, that
 \begin{equation}\label{eq:radial:nondeg}
  \langle u, y - x_0 \rangle \neq 0, ~ \forall ~ y \in K \setminus \{x_0\}, ~ \forall ~ u \in T(y),
 \end{equation}
 and that $T$ restricted to every ray $R_d$ ($d \in \mathbb{R}^n \setminus \{0\}$) is strongly quasimonotone with modulus $\gamma$. Then $T$ is strongly star monotone at $x_0$ with modulus $\gamma > 0$.
\end{theorem}

\begin{proof}
 Fix $y \in K \setminus \{x_0\}$ and let $d = y - x_0$. Since $K$ is star-shaped at $x_0$, the segment $[x_0, y]$ is contained in $K$. For $t \in \, ]0,1[$, let $x_t := x_0 + t(y - x_0)$, so that $x_0, x_t, y \in R_d$.

 Take any $u \in T(x_t)$, which exists since $T$ has nonempty values. We claim that $\langle u, x_0 - x_t \rangle < 0$. Indeed, if $\langle u, x_0 - x_t \rangle > 0$, then strong quasimonotonicity of $T$ on $R_d$ applied with $x = x_t$ and target point $x_0$ yields $\langle v, x_0 - x_t \rangle \geq \frac{\gamma}{2}\|x_0 - x_t\|^2$ for every $v \in T(x_0)$. Taking $v = 0 \in T(x_0)$ gives $0 \geq \frac{\gamma}{2} t^2 \|y - x_0\|^2 > 0$, a contradiction. On the other hand, $\langle u, x_0 - x_t \rangle \neq 0$ by \eqref{eq:radial:nondeg}, since $x_t \neq x_0$. Hence $\langle u, x_0 - x_t \rangle < 0$.

 Since $x_t - x_0 = t(y - x_0)$ and $y - x_t = (1-t)(y - x_0)$, we obtain
 \[
  \langle u, y - x_t \rangle = \frac{1-t}{t} \langle u, x_t - x_0 \rangle > 0.
 \]
 Then strong quasimonotonicity of $T$ on $R_d$ applied with $x = x_t$ and target point $y$ yields
 \[
  \langle v, y - x_t \rangle \geq \frac{\gamma}{2}\|y - x_t\|^2, ~ \forall ~ v \in T(y).
 \]
 Taking the limit as $t \to 0^+$ (so $x_t \to x_0$), we have
 \[
  \langle v, y - x_0 \rangle \geq \frac{\gamma}{2}\|y - x_0\|^2, ~ \forall ~ v \in T(y).
 \]
 Since $y \in K \setminus \{x_0\}$ was arbitrary, $T$ is strongly star monotone at $x_0$ with modulus $\gamma > 0$.
\end{proof}

\begin{remark} \label{no:relation:set}
 \begin{itemize}
  \item[$(i)$] Condition \eqref{eq:radial:nondeg} is necessary for strong star monotonicity. Indeed, if $T$ is strongly star monotone at $x_0$ with modulus $\gamma > 0$, then $\langle u, y - x_0 \rangle \geq \frac{\gamma}{2}\|y - x_0\|^2 > 0$ for every $y \in K \setminus \{x_0\}$ and every $u \in T(y)$. Furthermore, \eqref{eq:radial:nondeg} implies $0 \notin T(y)$ for every $y \neq x_0$, hence under the assumptions of Theorem \ref{theo:quasi:star:ray} we have $\zer(T) = \{x_0\}$.

  \item[$(ii)$] Condition \eqref{eq:radial:nondeg} cannot be dropped in Theorem \ref{theo:quasi:star:ray}, not even for single-valued continuous operators with a unique zero. Indeed, let $n = 2$, $x_0 = 0$, $\gamma > 0$, $J(y_1, y_2) := (-y_2, y_1)$, and
  \[
   T(y) := \frac{\gamma}{2} \left( 1 - \frac{1}{\|y\|} \right)^{+} y + Jy ~ \text{ if } y \neq 0, \qquad T(0) := 0.
  \]
  Then $T$ is continuous on $\mathbb{R}^2$ and $\zer(T) = \{0\}$, because $Jy \perp y$ and $Jy \neq 0$ for $y \neq 0$. For $\|d\| = 1$ and $x = ad$, $y = bd$ with $a, b \geq 0$, we have $\langle T(x), y - x \rangle = (b - a) \frac{\gamma}{2}(a - 1)^{+}$, which is positive only if $b > a > 1$, and in that case
  \[
   \langle T(y), y - x \rangle = (b-a) \frac{\gamma}{2} (b - 1) \geq \frac{\gamma}{2} (b - a)^2.
  \]
  Thus, $T$ restricted to every ray $R_d$ is strongly quasimonotone with modulus $\gamma$. However, $\langle T(y), y \rangle = 0$ for every $0 < \|y\| \leq 1$, so $T$ is not strongly star monotone at $0$. Of course, \eqref{eq:radial:nondeg} fails at every point $0 < \|y\| \leq 1$.

  \item[$(iii)$] When $n = 1$, condition \eqref{eq:radial:nondeg} simply means that $0 \notin T(y)$ for every $y \neq x_0$. Hence, in the one-dimensional case, Theorem \ref{theo:quasi:star:ray} holds under the sole assumptions $\zer(T) = \{x_0\}$ and strong quasimonotonicity of $T$ on both rays $[x_0, +\infty[ \, \cap K$ and $]-\infty, x_0] \cap K$.

  \item[$(iv)$] Theorem \ref{theo:quasi:star:ray} requires $\gamma > 0$. When $\gamma = 0$, quasimonotonicity and star monotonicity do not imply one another. For example, $T(x) = - \lvert x \rvert$ on $\mathbb{R}$ is quasimonotone but not star monotone, while $T = \nabla f$ from Remark \ref{no:rel:mon} (see also Remark \ref{no:rel:mon:robust} for a further instance drawn from robust statistics) is star monotone at $(0, 0)$ but not quasimonotone.
 \end{itemize}
\end{remark}

As a result, we have:

\begin{theorem} \label{prop:ray:equiv}
 Let $K \subset \mathbb{R}^n$ be star-shaped at $x_0 \in K$, $T: K \rightrightarrows \mathbb{R}^n$ be a set-valued operator with nonempty values such that $0 \in T(x_0)$, and $\gamma > 0$. Then the following are equivalent.
 \begin{enumerate}
  \item[$(a)$] $T$ is strongly star monotone at $x_0$ with modulus $\gamma$.

  \item[$(b)$] $T$ satisfies \eqref{eq:radial:nondeg}, and the restriction of $T$ to every ray $R_d$ ($d \neq 0$) is strongly quasimonotone with modulus $\gamma$.
 \end{enumerate}
\end{theorem}

\begin{proof}
 $(a) \Rightarrow (b)$ follows from Proposition \ref{prop:ray:necessity} and Remark \ref{no:relation:set}$(i)$, while $(b) \Rightarrow (a)$ is Theorem \ref{theo:quasi:star:ray}.
\end{proof}

\begin{remark} \label{rem:fails:strong}
 The assumption $0 \in T(x_0)$ cannot be dropped in Theorem \ref{prop:ray:equiv}. Indeed, let $K = \, ]-1/2, +\infty[$, $x_0 = 0$, $\gamma > 0$, and $T(x) = \gamma (x + 1)$. Then $T$ is strongly monotone with modulus $2\gamma$, hence strongly quasimonotone with modulus $\gamma$ on every ray, and \eqref{eq:radial:nondeg} holds since $T(x) > 0$ on $K$. However, $0 \notin T(0)$ and $\langle T(x), x \rangle < 0$ for every $x \in \, ]-1/2, 0[$, so $T$ is not star monotone at $0$.
\end{remark}

\section{Applications in Optimization and Variational Inequalities}\label{sec:04}

\subsection{The Separable Case}

In this subsection, we study separable star monotone operators (see \cite{CH1,CH2} for the quasimonotone case).

\begin{theorem} \label{extension:DK:set}
Let $C_1 \subset \mathbb{R}^{n_1}, \ldots, C_p \subset \mathbb{R}^{n_p}$ be sets, with $n_i \geq 1$ and $\sum_{i=1}^p n_i = n$. Let $T_i : C_i \rightrightarrows \mathbb{R}^{n_i}$ be set-valued operators for every $i = 1,\ldots,p$, and define $T : C_1 \times \cdots \times C_p \rightrightarrows \mathbb{R}^{n}$ by
\begin{equation}\label{sepa:oper}
 T(x_1,\ldots,x_p) := T_1(x_1) \times \cdots \times T_p(x_p).
\end{equation}
Let $x_0 = (x_1^0,\ldots,x_p^0) \in C_1 \times \cdots \times C_p$ and $\gamma \geq 0$. Then the following assertions hold.
\begin{itemize}
 \item[$(a)$] If every $T_i$ is (strongly) star monotone at $x_i^0$ with modulus $\gamma$, then $T$ is (strongly) star monotone at $x_0$ with modulus $\gamma$.

 \item[$(b)$] Suppose that $T_i(x_i^0) \neq \emptyset$ for every $i = 1, \ldots, p$. If $T$ is (strongly) star monotone at $x_0$ with modulus $\gamma$, then every $T_i$ is (strongly) star monotone at $x_i^0$ with modulus $\gamma$.
\end{itemize}
In particular, if $T_i(x_i^0) \neq \emptyset$ for every $i$ (e.g., if $0 \in T_i(x_i^0)$ for every $i$), then $T$ is (strongly) star monotone at $x_0$ with modulus $\gamma$ if and only if every $T_i$ is (strongly) star monotone at $x_i^0$ with modulus $\gamma$.
\end{theorem}

\begin{proof}
$(a)$: Suppose every $T_i$ is (strongly) star monotone at $x_i^0$ with modulus $\gamma \geq 0$. Let $y = (y_1,\ldots,y_p) \in C_1 \times \cdots \times C_p$ and $v = (v_1, \ldots, v_p) \in T(y)$, so $v_i \in T_i(y_i)$ for every $i$. Then,
\[
 \langle v, y-x_0 \rangle = \sum_{i=1}^p \langle v_i, y_i - x_i^0 \rangle \geq \sum_{i=1}^p \frac{\gamma}{2}\|y_i-x_i^0\|^2 = \frac{\gamma}{2}\|y-x_0\|^2.
\]
Since $y \in C_1\times\cdots\times C_p$ and $v \in T(y)$ were arbitrary, $T$ is (strongly) star monotone at $x_0$ with modulus $\gamma$.

$(b)$: Suppose $T$ is (strongly) star monotone at $x_0$ with modulus $\gamma \geq 0$. Fix $i \in \{1,\ldots,p\}$, $y_i \in C_i$ and $v_i \in T_i(y_i)$. Take $y_j := x_j^0$ for $j \neq i$, and let $y := (y_1,\ldots,y_p)$. Since $T_j(x_j^0) \neq \emptyset$, we can pick $v_j \in T_j(x_j^0)$ for every $j \neq i$. Then $v := (v_1,\ldots,v_p) \in T(y)$, and by the star monotonicity of $T$ at $x_0$, we have
\[
 \langle v_i, y_i - x_i^0 \rangle = \sum_{k=1}^p \langle v_k, y_k - x_k^0 \rangle = \langle v, y - x_0 \rangle \geq \frac{\gamma}{2}\|y-x_0\|^2 = \frac{\gamma}{2}\|y_i - x_i^0\|^2.
\]
As $y_i \in C_i$ and $v_i \in T_i(y_i)$ were arbitrary, $T_i$ is (strongly) star monotone at $x_i^0$ with modulus $\gamma$.
\end{proof}

\begin{remark}\label{rem:nonempty:set}
 The assumption $T_i(x_i^0) \neq \emptyset$ in Theorem \ref{extension:DK:set}$(b)$ cannot be dropped. Indeed, let $p = 2$, $C_1 = C_2 = \mathbb{R}$, $x_0 = (0,0)$, $\gamma = 0$, and
 \[
  T_1(y_1) := -\frac{y_1}{1 + y_1^2}, \qquad
  T_2(y_2) :=
  \begin{cases}
   \{2/y_2\} & {\rm if } ~~ y_2 \neq 0, \\
   \emptyset & {\rm if } ~~ y_2 = 0.
  \end{cases}
 \]
 If $y_2 = 0$, then $T(y) = \emptyset$ and \eqref{eq:star:mon} holds trivially, while if $y_2 \neq 0$, then for $v \in T(y)$ we have $\langle v, y \rangle = 2 - \frac{y_1^2}{1 + y_1^2} \geq 1 > 0$. Hence, $T$ is star monotone at $(0,0)$. However, $\langle T_1(y_1), y_1 \rangle = -\frac{ y_1^2}{ 1+y_1^2} < 0$ for every $y_1 \neq 0$, so $T_1$ is not star monotone at $0$.
\end{remark}

As a consequence, we present a partial extension of \cite[Theorem 3]{Lara-10}, related to the Debreu--Koopmans problem (see \cite{DK}).

\begin{corollary} \label{coro:01:set}
 Let $T$ be a separable set-valued operator as in Theorem \ref{extension:DK:set}, where each $C_i \subset \mathbb{R}^{n_i}$ is star-shaped at $x_i^0$, and let $\gamma > 0$. Suppose that, for every $i = 1, \ldots, p$, $T_i$ has nonempty values, $0 \in T_i(x_i^0)$,
 \begin{equation}\label{eq:radial:nondeg:i}
  \langle u_i, y_i - x_i^0 \rangle \neq 0, ~ \forall ~ y_i \in C_i \setminus \{x_i^0\}, ~ \forall ~ u_i \in T_i(y_i),
 \end{equation}
 and $T_i$ restricted to every ray $\{x_i^0 + s d_i : s \geq 0\} \cap C_i$ ($d_i \in \mathbb{R}^{n_i} \setminus \{0\}$) is strongly quasimonotone with modulus $\gamma$. Then the following assertions hold.
 \begin{itemize}
  \item[$(a)$] $T$ is strongly star monotone at $x_0 = (x_1^0, \ldots, x_p^0)$ with modulus $\gamma$.

  \item[$(b)$] If $n_i = 1$ for every $i$, then the conclusion of $(a)$ holds with \eqref{eq:radial:nondeg:i} replaced by $\zer(T_i) = \{x_i^0\}$ for every $i = 1, \ldots, p$.
 \end{itemize}
\end{corollary}

\begin{proof}
 $(a)$: Fix $i \in \{1, \ldots, p\}$. Since $C_i$ is star-shaped at $x_i^0$, $T_i$ has nonempty values, $0 \in T_i(x_i^0)$, $T_i$ satisfies \eqref{eq:radial:nondeg:i}, and $T_i$ is strongly quasimonotone with modulus $\gamma > 0$ along rays from $x_i^0$, Theorem \ref{theo:quasi:star:ray} yields that $T_i$ is strongly star monotone at $x_i^0$ with modulus $\gamma$. Since this holds for every $i = 1, \ldots, p$, the implication $(\Leftarrow)$ in Theorem \ref{extension:DK:set}$(a)$ yields that $T$ is strongly star monotone at $x_0$ with modulus $\gamma$.

 $(b)$: If $n_i = 1$, then for every $y_i \in C_i \setminus \{x_i^0\}$ and $u_i \in T_i(y_i)$ we have $\langle u_i, y_i - x_i^0 \rangle = u_i (y_i - x_i^0)$, which is nonzero if and only if $u_i \neq 0$. Hence, \eqref{eq:radial:nondeg:i} is equivalent to $\zer(T_i) \subset \{x_i^0\}$ (see Remark \ref{no:relation:set}$(iii)$), and the result follows from $(a)$.
\end{proof}

\begin{remark}
 \begin{itemize}
  \item[$(i)$] In the quasimonotone case (i.e., when $\gamma = 0$), the implication in Corollary \ref{coro:01:set} does not hold. Indeed, even if every $T_i$ is quasimonotone and has a unique zero, the separable operator $T$ need not be star monotone. Take $p = 2$ and $T_1(x) = T_2(x) = -|x|$ for $x \in \mathbb{R}$. Every $T_i$ is continuous on $\mathbb{R}$, quasimonotone, and has a unique zero at $0$, but it is not star monotone at any point (see Remark \ref{no:relation:set}$(iv)$). Define $T: \mathbb{R}^2 \to \mathbb{R}^2$ by
  $$T(x_1,x_2) = (T_1(x_1), T_2(x_2)) = (-|x_1|, -|x_2|).$$
  By Theorem \ref{extension:DK:set}, $T$ is star monotone at a point $x_0 = (x_1^0, x_2^0)$ if and only if each $T_i$ is star monotone at $x_i^0$. Since neither component is star monotone at any point, $T$ is not star monotone at any point.

  \item[$(ii)$] Theorem \ref{extension:DK:set} and Corollary \ref{coro:01:set} are the star monotonicity analogues of the classical separability results in \cite{CH1,CH2}, which characterize when a separable product $T = (T_1, \ldots, T_p)$ is quasimonotone via their monotonicity indices $m(F_i)$, in the spirit of the Debreu--Koopmans theorem \cite{DK} for separable sums of quasiconvex functions.

  Since star monotonicity is an anchored, one-point condition, it decouples exactly across factors, with no exceptional-factor case and no index arithmetic. In this sense, Theorem \ref{extension:DK:set} is closer in spirit to the separable characterization of star quasiconvexity given in \cite[Theorem 3]{Lara-10} and to the resolution of the Debreu--Koopmans problem in \cite[Corollary 4]{Lara-10} (see also Remark \ref{rem:important}$(ii)$), now lifted to the operator setting. In particular, Corollary \ref{coro:01:set}$(b)$ shows that, in the scalar case $n_i = 1$, ray-wise strong quasimonotonicity and a unique zero in each factor are enough to guarantee strong star monotonicity of the separable operator.
 \end{itemize}
\end{remark}

\subsection{Generalized Convexity and Generalized Monotonicity}

%We characterize nonsmooth (strongly) star quasiconvexity functions via their Clarke subdifferential \cite{clarke}. To that end, 
Let us recall the following notions from nonsmooth analysis.

Let $f:\R^n \to \R$ be locally Lipschitz around $x \in \R^{n}$. The Clarke directional derivative of $f$ at $x$ in direction $d \in \mathbb{R}^n$ is defined by
\[
 f^{\circ} (x; d) := \mathop{\limsup}\limits_{ \substack{y \to x \\ t\downarrow 0}} \frac{f(y+td) - f(y)}{t},
\]
 and the Clarke subdifferential of $f$ at $x$ is given by
 \begin{equation*}%\label{def:clarke} 
  \partial^{c} f(x) := \{v \in \mathbb R^n: \, \langle v, d \rangle \le f^\circ(x; d), ~ \forall ~ d \in \mathbb R^n\},
 \end{equation*} 
 Moreover, by \cite[Proposition 2.1.2]{clarke},
\begin{equation}\label{prop:direc}
 f^{\circ} (x; d) = \max_{\xi \in \partial^{c} f(x)} \langle \xi, d \rangle, ~ \forall ~ d \in \mathbb R^{n}.
\end{equation}
% A point $x \in \R^{n}$ is called a Clarke critical point of $f$ if $0 \in \partial^{c} f(x)$.

Let us also recall the following results.

\begin{lemma} {\rm (Lebourg's Mean Value Theorem \cite[Theorem 2.3.7]{clarke})} \label{lem:mvt}
 Let $f$ be locally Lipschitz continuous on an open convex set $K \subset \mathbb{R}^n$. For every $x, y \in K$, there exists $u \in \, ]x, y[$ and $\xi \in \partial^c f(u)$ such that
 $$f(x) - f(y) = \langle \xi, x - y\rangle.$$
\end{lemma}

\begin{lemma} {\rm (\cite[Theorem 2.5.1]{clarke})} \label{lem:clarke}
 Let $f: K \to \mathbb{R}$ be locally Lipschitz continuous on an open convex set $K$ and $D \subset K$ the set of points where $f$ is differentiable. Then for every $y \in K$,
 \[
  \partial^c f(y) = {\rm co}(S(y)),
 \]
 where $S(y) := \left\{ \xi^* \in \mathbb{R}^{n}: \exists\, \{y_k\}_{k} \subset D,\ y_k \to y,\ \nabla f(y_k) \to \xi^{*} \right\}$.
\end{lemma}

In the following result, we characterize (strongly) star quasiconvexity of a function via the strong star monotonicity of its Clarke subdifferential. Other characterizations for generalized convexity notions using Clarke subdifferential may be found in \cite{FLG}.

\begin{theorem} \label{thm:genconv}
 Let $f: \mathbb{R}^n \to \mathbb{R}$ be locally Lipschitz continuous on an open convex set $K \subset \mathbb{R}^{n}$, and suppose $\overline{x} \in {\rm argmin}_{K}\,f$. Then $f$ is (strongly) star qua\-siconvex on $K$ with respect to $\overline{x}$ with modulus $\gamma \geq 0$ if and only if 
\begin{equation}
 \langle \xi_{y}, y - \overline{x} \rangle \geq \frac{\gamma}{2} \|y - \overline{x}\|^2, ~ \forall ~ y \in K, ~ \forall ~ \xi_{y} \in \partial^c f(y),
\end{equation}
i.e., the Clarke subdifferential $\partial^c f$ is (strongly) star monotone at $\overline{x}$ with modulus $\gamma \geq 0$ on $K$.
\end{theorem}

\begin{proof}
 ($\Rightarrow$) Let $z \in K \cap D$ be any point where $f$ is differentiable. By Definition \ref{def:starqc}, we have for every $\lambda \in \, ]0, 1[$ that
\begin{align*}
 & \hspace{2.2cm} f(z + \lambda(\overline{x} - z)) \leq f(z) - \lambda(1-\lambda)\frac{\gamma}{2}\|z - \overline{x}\|^2 \\
 & \iff \frac{f(z + \lambda ( \overline{x} - z)) - f(z)}{\lambda} \leq -(1-\lambda) \frac{\gamma}{2} \|z - \overline{x}\|^2.
\end{align*}
Taking the limit as $\lambda \downarrow 0$ and using the differentiability of $f$ at $z$,
\begin{equation}\label{tag:2}
 \langle \nabla f(z), \overline{x} - z \rangle \leq -\frac{\gamma}{2}\|z - \overline{x}\|^2 \iff \langle \nabla f(z), z - \overline{x} \rangle \geq \frac{\gamma}{2}\|z - \overline{x}\|^2. 
\end{equation}

Now, let $y \in K$ be arbitrary and $\xi \in \partial^c f(y)$. By Lemma \ref{lem:clarke}, we have $\partial^c f(y) = {\rm co}(S(y))$. For any generator $\xi^* \in S(y)$, there exists a sequence $\{y_k\}_{k} \subset D$ with $y_k \to y$ such that $\nabla f(y_k) \to \xi^*$. Evaluating \eqref{tag:2} at each $y_k$, 
\[
 \langle \nabla f(y_k), y_k - \overline{x} \rangle \geq \frac{\gamma}{2}\|y_k - \overline{x}\|^2.
\]
Taking $k \to \infty$ yields $\langle \xi^*, y - \overline{x} \rangle \geq \frac{\gamma}{2}\|y - \overline{x}\|^2$. Since $\xi \in \partial^c f(y)$ is a convex combination of such generators, $\xi = \sum_{i=1}^m \lambda_i \xi_i^*$ with $\lambda_i \geq 0$ and $\sum_{i=1}^m \lambda_i = 1$, linearity of the inner product gives
\[
 \langle \xi, y - \overline{x} \rangle = \sum_{i=1}^m \lambda_i \langle \xi_i^*, y - \overline{x} \rangle \geq \sum_{i=1}^m \lambda_i \frac{\gamma}{2}\|y - \overline{x}\|^2 = \frac{\gamma}{2}\|y - \overline{x}\|^2.
\]
Hence, $\partial^c f$ is (strongly) star monotone with modulus $\gamma \geq 0$.

\medskip
\noindent ($\Leftarrow$) Assume that $\partial^c f$ is (strongly) star monotone with modulus $\gamma \geq 0$. Fix $y \in K$ and $\lambda \in \, ]0, 1]$ and set $y_\lambda := \lambda \overline{x} + (1-\lambda)y$. By Lebourg's Mean Value Theorem (Lemma \ref{lem:mvt}), there exists $h \in \, ]0, \lambda[$ and $\xi_h \in \partial^c f(y_h)$, where $y_h = h\overline{x} + (1-h)y$, such that
\begin{equation}\label{tag:3}
 f(y) - f(y_\lambda) = \langle \xi_h, y - y_\lambda \rangle = \lambda \langle \xi_h, y - \overline{x} \rangle. 
\end{equation}

Since $y_h - \overline{x} = (1-h)(y - \overline{x})$, by the strong star monotonicity of $\partial^c f$ at $\overline{x}$,
$$\langle \xi_h, y_h - \overline{x} \rangle \geq \frac{\gamma}{2} \|y_h - \overline{x}\|^2 ~ \Longrightarrow ~ \langle \xi_h, (1-h)(y - \overline{x}) \rangle \geq \frac{\gamma}{2} (1-h)^2 \|y - \overline{x}\|^2.$$
Dividing by $1-h$ yields
\begin{equation}\label{tag:5}
 \langle \xi_h, y - \overline{x} \rangle \geq \frac{\gamma}{2} (1-h) \|y - \overline{x}\|^2. 
\end{equation}
Since $h < \lambda$, we have $1-h > 1-\lambda$. Substituting \eqref{tag:5} into \eqref{tag:3} gives
$$f(y) - f(y_\lambda) = \lambda \langle \xi_h, y - \overline{x} \rangle \geq \lambda (1-\lambda) \frac{\gamma}{2} \|y - \overline{x}\|^2,$$ 
which rearranges to $f(\lambda \overline{x} + (1-\lambda)y) \leq f(y) - \lambda(1-\lambda) \frac{ \gamma}{2}\|y - \overline{x}\|^2$. Thus $f$ is (strongly) star quasiconvex with modulus $\gamma \geq 0$.
\end{proof}

\subsection{Applications to Minty and Stampacchia Variational Inequalities}

We apply the star monotonicity framework to establish existence and equivalence results for Minty and Stampacchia Variational Inequalities. 

Given a set-valued operator $T: K \rightrightarrows \mathbb{R}^n$ on a set $K \subset \mathbb{R}^n$, we recall the classical variational inequality problems, together with a strong version of the Minty Variational Inequality.
\begin{itemize}
 \item[$(i)$] Stampacchia Variational Inequality (SVI): Find $x_0 \in K$ such that there exists $u_0 \in T(x_0)$ satisfying
 \[
  \langle u_0, y - x_0 \rangle \geq 0, ~ \forall ~ y \in K.
  \]
 \item[$(ii)$] Minty Variational Inequality (MVI): Find $x_0 \in K$ such that
  \[
  \langle v, y - x_0 \rangle \geq 0, ~ \forall ~ y \in K, ~ \forall ~ v \in T(y).
  \]
 \item[$(iii)$] Strong Minty Variational Inequality with modulus $\gamma > 0$ (Strong MVI): Find $x_0 \in K$ such that
 \[
  \langle v, y - x_0 \rangle \geq \frac{\gamma}{2}\|y - x_0\|^2, ~ \forall ~ y \in K, ~ \forall ~ v \in T(y).
 \]
 \end{itemize}

Note that $x_0$ solves the Strong MVI with modulus $\gamma > 0$ precisely when $T$ is strongly star monotone at $x_0$ with modulus $\gamma$. Likewise, $x_0$ solves the MVI precisely when $T$ is star monotone at $x_0$ with modulus $\gamma = 0$.

First, we establish the complete equivalence between zero points, solutions to SVI, solutions to MVI, and strong star monotonicity in the setting where $\gamma > 0$.

\begin{theorem}[Minty--Stampacchia Equivalence for $\gamma > 0$] \label{theo:minty:strong}
 Let $K \subset \mathbb{R}^n$ be an open set star-shaped at $x_0 \in K$, $\gamma > 0$, and $T: K \rightrightarrows \mathbb{R}^n$ be a set-valued operator with nonempty values which is weakly hemicontinuous at $x_0$ and satisfies the radial nondegeneracy condition \eqref{eq:radial:nondeg}. Suppose that $T$ restricted to every ray $R_d := \{x_0 + sd : s \ge 0\} \cap K$ ($d \neq 0$) is strongly quasimonotone with modulus $\gamma$. Then the following assertions are equivalent.
 \begin{enumerate}
  \item[$(a)$] $T$ is strongly star monotone at $x_0$ with modulus $\gamma$.

  \item[$(b)$] $x_0$ solves the Strong Minty Variational Inequality with modulus $\gamma$.

  \item[$(c)$] $x_0$ solves the Stampacchia Variational Inequality (SVI).

  \item[$(d)$] $x_0 \in \zer(T)$.
 \end{enumerate}
 In this case, $\zer(T) = \{x_0\}$.
\end{theorem}

\begin{proof}
 $(a) \Leftrightarrow (b)$: This is the definition of strong star monotonicity at $x_0$.

 $(b) \Rightarrow (d)$: If $x_0$ solves the Strong MVI, then $T$ is strongly star monotone at $x_0$. Since $T$ is weakly hemicontinuous at $x_0$ and $T(x_0) \neq \emptyset$, Theorem \ref{theo:zero}$(a)$ guarantees that $0 \in T(x_0)$.

 $(d) \Rightarrow (c)$: Take $u_0 = 0 \in T(x_0)$. Then $\langle u_0, y - x_0 \rangle = 0$ for all $y \in K$, so $x_0$ solves SVI.

 $(c) \Rightarrow (d)$: Suppose $x_0$ solves SVI, so there exists $u_0 \in T(x_0)$ such that $\langle u_0, y - x_0 \rangle \geq 0$ for all $y \in K$. Fix any direction $d \in \mathbb{R}^n \setminus \{0\}$. Since $K$ is open, there exists $t > 0$ such that $x_0 \pm t d \in K$. Taking $y = x_0 + td$ and $y = x_0 - td$ yields $\langle u_0, d \rangle \geq 0$ and $\langle u_0, d \rangle \leq 0$, respectively. Thus, $\langle u_0, d \rangle = 0$ for all $d \in \mathbb{R}^n$, which implies $u_0 = 0 \in T(x_0)$.

 $(d) \Rightarrow (a)$: Since $K$ is star-shaped at $x_0$, $T$ has nonempty values, $0 \in T(x_0)$, $T$ satisfies \eqref{eq:radial:nondeg}, and $T$ restricted to every ray $R_d$ is strongly quasimonotone with modulus $\gamma > 0$, Theorem \ref{theo:quasi:star:ray} yields that $T$ is strongly star monotone at $x_0$ with modulus $\gamma$.

 Finally, \eqref{eq:radial:nondeg} implies $0 \notin T(y)$ for every $y \in K \setminus \{x_0\}$, hence $\zer(T) = \{x_0\}$ whenever $(d)$ holds.
\end{proof}

\begin{remark} \label{rem:minty:hyp}
 \begin{itemize}
  \item[$(i)$] The hypotheses of Theorem \ref{theo:minty:strong} do not imply any of the equivalent assertions $(a)$--$(d)$. Indeed, the operator $T(x) = \gamma(x+1)$ on $K = \, ]-1/2, +\infty[$ from Remark \ref{rem:fails:strong} is continuous, satisfies \eqref{eq:radial:nondeg} at $x_0 = 0$, and is strongly quasimonotone with modulus $\gamma$ on every ray, but $0 \notin T(0)$. Hence, Theorem \ref{theo:minty:strong} is a genuine characterization: under ray-wise strong quasimonotonicity and radial nondegeneracy, strong star monotonicity at $x_0$, the Strong MVI, the SVI, and the zero condition at $x_0$ are mutually equivalent.

  \item[$(ii)$] The radial nondegeneracy condition \eqref{eq:radial:nondeg} cannot be dropped in Theorem \ref{theo:minty:strong}. Indeed, the operator $T$ in Remark \ref{no:relation:set}$(ii)$ is continuous on $\mathbb{R}^2$ and strongly quasimonotone with modulus $\gamma$ on every ray from $0$, and satisfies $0 \in T(0)$, so $(c)$ and $(d)$ hold. However, $T$ is not strongly star monotone at $0$, so $(a)$ and $(b)$ fail.
 \end{itemize}
\end{remark}

We close this subsection with the $\gamma = 0$ counterpart of Theorem \ref{theo:minty:strong}, which follows directly from results already established and therefore requires no separate proof.

\begin{remark} \label{rem:minty:zero}
 The MVI is, by definition, exactly the case $\gamma = 0$ of the Strong MVI, i.e., $x_0$ solves the MVI precisely when $T$ is star monotone at $x_0$ with modulus $\gamma = 0$. Consequently, the $\gamma > 0$ equivalence of Theorem \ref{theo:minty:strong} degenerates, at $\gamma = 0$, to two implications that are already available. On the one hand, if $0 \in T(x_0)$ and $T$ is pseudomonotone on $K$, then $x_0$ solves the MVI, by Proposition \ref{pseudo:star:set} with $\gamma = 0$. On the other hand, if $K$ is open, $T(x_0) \neq \emptyset$, $T$ is weakly hemicontinuous at $x_0$, and $x_0$ solves the MVI, then $0 \in T(x_0)$, by Theorem \ref{theo:zero}$(a)$. Furthermore, if $K$ is open, then $x_0 \in \zer(T)$ if and only if $x_0$ solves the SVI, without any monotonicity or continuity assumption on $T$: indeed, steps $(d) \Rightarrow (c)$ and $(c) \Rightarrow (d)$ in the proof of Theorem \ref{theo:minty:strong} use only the openness of $K$.
\end{remark}

\subsection{Applications to First-Order Methods}

As a first application of strong star monotonicity to algorithmic convergence, we recall the following linear-convergence result for a general fixed-point scheme.

\begin{proposition}[Linear Convergence of a Fixed-Point Scheme]\label{prop:linear:conv}
 Let $T: \mathbb{R}^n \to \mathbb{R}^n$ be single-valued, strongly star monotone at $x_0$ with modulus $\gamma > 0$ and $T(x_{0}) = 0$, and suppose there exists $L > 0$ such that $\|T(x)\| \leq L\|x - x_0\|$ for all $x \in \mathbb{R}^n$ (note $\gamma \leq 2L$ automatically, by Cauchy--Schwarz and star monotonicity). Fix $0 < \eta < \gamma/L^2$, take $x^0 \in \mathbb{R}^n$, and define $x^{k+1} := x^k - \eta\, T(x^k)$. Then,
 \begin{equation}\label{linear:conv}
  \|x^{k+1} - x_0\|^2 \leq \left(1 - \eta\gamma + \eta^2 L^2\right) \|x^k - x_0\|^2, ~ \forall ~ k \in \mathbb{N},
 \end{equation}
 so $x^k \to x_0$ linearly. The choice $\eta = \gamma/(2L^2)$ minimizes the contraction factor, giving rate $1 - \gamma^2/(4L^2)$.
\end{proposition}

\begin{proof}
 Expand directly
 $$\|x^{k+1}-x_0\|^2 = \|x^k - x_0\|^2 - 2\eta\langle T(x^k), x^k - x_0\rangle + \eta^2 \|T(x^k)\|^2.$$
 By the strong star monotonicity, $\langle T(x^k), x^k - x_0\rangle \geq \frac{\gamma}{2}\|x^k-x_0\|^2$, while by the growth bound we have $\|T(x^k)\|^2 \leq L^2\|x^k-x_0\|^2$. Hence, the result follows by substituting both in the previous relation.

 Furthermore, minimizing $1 - \eta \gamma + \eta^2 L^2$ over $\eta > 0$ gives $\eta = \gamma/(2L^2)$, which completes the proof.
\end{proof}

While Proposition \ref{prop:linear:conv} applies to single-valued operators satisfying a global Lipschitz-type growth bound, it does not cover the nonsmooth setting of Section 4.2, where $T = \partial^c f$ is generally set-valued and only locally (not necessarily globally, from $x_0$) controlled. We close this section by combining Proposition \ref{prop:linear:conv}'s technique with Theorem \ref{thm:genconv} to analyze the classical subgradient method for minimizing a nonsmooth star quasiconvex function.

\begin{theorem}[Linear Convergence to a Neighborhood via the Subgradient Method]\label{theo:subgrad}
 Let $f: \mathbb{R}^n \to \mathbb{R}$ be locally Lipschitz continuous and strongly star quasiconvex with respect to $\overline{x} \in \operatorname{argmin}_{\mathbb{R}^n} f$ with modulus $\gamma > 0$. Let $x^0 \in \mathbb{R}^n \setminus \{\overline{x}\}$, $R := \|x^0 - \overline{x}\|$, and
 \[
  L := \sup \left\{ \|\xi\| : \xi \in \partial^c f(y), ~ y \in \overline{\mathbb{B}}(\overline{x}, R) \right\},
 \]
 which is finite since $\partial^c f$ is locally bounded {\rm \cite[Proposition 2.1.2]{clarke}}, and positive since $L \geq \frac{\gamma}{2} R$ by Theorem \ref{thm:genconv} and Proposition \ref{prop:01:set}$(d)$. Fix $0 < \eta < 1/\gamma$ with $\eta \leq \gamma R^2/L^2$, and generate the sequence
 \[
  x^{k+1} := x^k - \eta\, \xi_k, \qquad \xi_k \in \partial^c f(x^k) \ \text{(any selection)}.
 \]
 Then $\{x^k\}_k \subset \overline{\mathbb{B}}(\overline{x}, R)$ and
 \begin{equation}\label{lin:clarke}
  \|x^k - \overline{x}\|^2 \leq (1-\eta \gamma)^k \|x^0 - \overline{x}\|^2 + \frac{\eta L^2}{\gamma} \left( 1 - (1-\eta\gamma)^k \right), ~ \forall ~ k \in \mathbb{N},
 \end{equation}
 so $\{x^k\}_{k}$ converges linearly to the ball $\overline{\mathbb{B}} (\overline{x}, \sqrt{\eta L^2/\gamma})$, i.e.
 \begin{equation}\label{ball:conv}
  \limsup_{k\to\infty} \|x^k - \overline{x}\|^2 \leq \frac{\eta L^2}{\gamma}.
 \end{equation}
\end{theorem}

\begin{proof}
 By Theorem \ref{thm:genconv}, $\partial^c f$ is strongly star monotone at $\overline{x}$ with modulus $\gamma > 0$. Suppose that $\|x^k - \overline{x}\| \leq R$ for some $k$. Then $\|\xi_k\| \leq L$ and $\langle \xi_k, x^k - \overline{x} \rangle \geq \frac{\gamma}{2}\|x^k - \overline{x}\|^2$, so
 \begin{align}
  \|x^{k+1} - \overline{x}\|^2 & = \|x^k - \overline{x} \|^2 - 2\eta \langle \xi_k, x^k - \overline{x} \rangle + \eta^2 \|\xi_k\|^2 \notag \\
  & \leq (1-\eta\gamma)\|x^k-\overline{x}\|^2 + \eta^2 L^2. \label{rec:subgrad}
 \end{align}
 In particular, $\|x^{k+1} - \overline{x}\|^2 \leq (1 - \eta\gamma) R^2 + \eta^2 L^2 \leq R^2$, since $\eta L^2 \leq \gamma R^2$. As $\|x^0 - \overline{x}\| = R$, induction shows that $\{x^k\}_k \subset \overline{\mathbb{B}}(\overline{x}, R)$ and that \eqref{rec:subgrad} holds for every $k \in \mathbb{N}$.

 Since $0 < \eta < 1/\gamma$, iterating \eqref{rec:subgrad} from $k=0$ gives, by summing the resulting geometric series,
 \begin{align*}
  \|x^k-\overline{x}\|^2 & \leq (1-\eta \gamma)^k \|x^0 - \overline{x}\|^2 + \eta^2L^2 \sum_{j=0}^{k-1}(1- \eta \gamma)^j \\
  & = (1-\eta \gamma)^k\|x^0-\overline{x} \|^2 + \frac{ \eta L^2}{\gamma} \left( 1- (1- \eta \gamma )^k \right).
 \end{align*}
 Since $0 < 1- \eta \gamma < 1$, the result follows by taking $\limsup_{k \to \infty}$.
\end{proof}

\begin{remark}
 Unlike Proposition \ref{prop:linear:conv}, exact linear convergence to $\overline{x}$ does not hold here with a constant step size: the growth bound $\|T(x)\|\leq L\|x-x_0\|$ used there forces $T$ to vanish at the target point, whereas a Clarke subgradient $\xi \in \partial^c f(\overline{x})$ need not be zero when $f$ has a kink at $\overline{x}$, only $0\in\partial^c f(\overline{x})$ is guaranteed. Theorem \ref{theo:subgrad} gives the standard trade-off of nonsmooth subgradient methods: a larger step $\eta$ yields faster initial linear decay but converges to a larger neighborhood of $\overline{x}$.
\end{remark}

In particular, we obtain the following exact convergence result.

\begin{corollary} \label{cor:subgrad:diminish}
 Let $f$, $\overline{x}$, $x^0$, $R$ and $L$ be as in Theorem \ref{theo:subgrad}, and let $\{\eta_k\}_{k} \subset \, ]0, 1/\gamma[$ satisfy $\eta_k \leq \gamma R^2/L^2$ for every $k$, $\sum_{k=0}^{\infty} \eta_k = + \infty$ and $\eta_k \to 0$. Generate the sequence
 $$x^{k+1} := x^k - \eta_k \xi_k, ~~~ \xi_k \in \partial^c f(x^k).$$
 Then, $x^k \to \overline{x}$.
\end{corollary}

\begin{proof}
 Arguing as in the proof of Theorem \ref{theo:subgrad} with $\eta_k$ instead of $\eta$, we obtain $\{x^k\}_k \subset \overline{\mathbb{B}}(\overline{x}, R)$ and, for $a_k := \|x^k - \overline{x}\|^2$,
 \[
  a_{k+1} \leq (1 - \eta_k \gamma) a_k + \eta_k^2 L^2, ~ \forall ~ k \in \mathbb{N}.
 \]
 Fix $\varepsilon > 0$. Since $\eta_k \to 0$, there exists $k_0$ such that $\eta_k L^2 \leq \gamma \varepsilon$ for all $k \geq k_0$. For such $k$,
 \[
  a_{k+1} - \varepsilon \leq (1 - \eta_k\gamma)(a_k - \varepsilon) + \eta_k (\eta_k L^2 - \gamma\varepsilon) \leq (1 - \eta_k\gamma)(a_k - \varepsilon),
 \]
 hence $(a_{k+1} - \varepsilon)^+ \leq (1 - \eta_k\gamma)(a_k - \varepsilon)^+$. Therefore, for every $k > k_0$,
 \[
  (a_k - \varepsilon)^+ \leq \prod_{j=k_0}^{k-1} (1 - \eta_j\gamma) \, (a_{k_0} - \varepsilon)^+ \leq \exp\Big( -\gamma \sum_{j=k_0}^{k-1} \eta_j \Big) (a_{k_0} - \varepsilon)^+ \to 0,
 \]
 since $\sum_j \eta_j = +\infty$. Thus $\limsup_{k} a_k \leq \varepsilon$, and since $\varepsilon > 0$ was arbitrary, $x^k \to \overline{x}$.
\end{proof}

\section{Conclusions}\label{sec:05}

We introduced star monotone operators and developed their basic theory, including calculus rules, existence and uniqueness of zeros, a global metric subregularity estimate, and a precise structural relationship to the classical monotone, pseudomonotone, and quasimonotone hierarchy. We further demonstrated the applicability of the notion through a Clarke subdifferential characterization of nonsmooth star quasiconvexity, a complete Minty--Stampacchia equivalence, linear convergence guarantees for first-order methods, and a separable characterization related to the Debreu--Koopmans problem.

%A natural direction for future work is to explore extensions of the classical Zarantonello--Minty surjectivity theory \cite{Zaran,Minty1962} to star monotone operators, and studying the resulting resolvent (its single-valuedness, nonexpansiveness, and continuity properties) would place star monotone operators within the classical operator-theoretic framework underlying the proximal point algorithm \cite{rock-1976}, and open the way to resolvent-based algorithms beyond the first-order methods considered here.

\section{Declarations}

%\subsection{Ethical Approval and Consent to participate}

%Not Applicable.

%\subsection{Consent for publication}

%Not Applicable.

%\subsection{Human and Animal Ethics}

%Not Applicable.

\subsection{Availability of Supporting Data}

No data sets were generated during the current study. 

\subsection{Author Contributions}

 All authors contributed equally to the study conception, design and implementation and wrote and corrected the manuscript.

\subsection{Conflicts of Interest}

There are no conflicts of interest or competing interests related to this manuscript.

\subsection{Funding}

This work has been partially su\-ppor\-ted by NAFOSTED---Vietnam through grant 101.01-2025.61 (Khanh), and by ANID---Chile through Fondecyt Regular 1241040 (Lara).

\subsection{Acknowledgments}

Parts of this work were carried out when F. Lara was vi\-si\-ting the Ton Duc Thang University (TDTU) in Ho Chi Minh City, Vietnam, during July 2026, and when both authors had a research stay at the Vietnam Institute for Advanced Study in Mathematics (VIASM), in Hanoi, Vietnam, during August-September 2026.
The authors wish to thank the above institutions for their hospitality.

\end{document}